\documentclass{article}
\usepackage[utf8]{inputenc}
\usepackage{multirow}
\usepackage{longtable,verbatim}
\usepackage{enumitem,  rotfloat, supertabular}
\usepackage{booktabs}
\usepackage{amsmath,hhline,dcolumn,layout,booktabs, pdflscape, amsthm}
\allowdisplaybreaks
\usepackage{amsfonts}
\usepackage[round]{natbib}

\floatstyle{ruled}

\usepackage{hyperref}
    \hypersetup{
        colorlinks   = true,
        citecolor    = blue,
        linkcolor=blue
    }
\usepackage{geometry}
\usepackage{adjustbox}
\usepackage{setspace}
\usepackage{multicol}
\usepackage{lipsum}
\usepackage{mwe}
\usepackage{style}
\usepackage[utf8]{inputenc}
\usepackage[most]{tcolorbox}
\usepackage{tikz}
\usepackage[utf8]{inputenc}
\usepackage[utf8]{inputenc}
\usepackage[utf8]{inputenc}
\usepackage{multirow}
\usepackage{longtable,verbatim}
\usepackage{enumitem,  rotfloat, supertabular}
\usepackage{booktabs}
\usepackage{amsmath,hhline,dcolumn,layout,booktabs, pdflscape, amsthm}
\usepackage{amsfonts}
\usepackage[round]{natbib}
\usepackage{mathtools}
\algrenewcommand\algorithmicrequire{\textbf{Input:}}
\algrenewcommand\algorithmicensure{\textbf{Output:}}
\algrenewcommand{\algorithmiccomment}[1]{\hskip3em$\rightarrow$ #1}

\floatstyle{ruled}
\floatname{algorithm}{Algorithm}
\renewcommand{\algorithmiccomment}[1]{// #1}
\usepackage{geometry}
\usepackage{hyperref}
    \hypersetup{
        colorlinks   = true,
        citecolor    = blue,
        linkcolor=blue  
    }
\tcbset{
    enhanced,
    boxrule=4pt, 
    colback=white,
    colframe=purple,
    halign=center,
    valign=center, 
    halign lower=center,
    valign lower =center,
    sharp corners=all,
    center title,
    lower separated=false,
    fonttitle=\sffamily\bfseries\large,
    }
\usepackage{adjustbox}
\usepackage{amsmath}

\algrenewcommand\algorithmicforall{\textbf{foreach}}

\usepackage{caption}
\usepackage{subcaption}
\usepackage{placeins}

\title{Joint Facility and Demand Location Problem}

\author{Ali Kaan Kurbanzade $^{a,b}$ $\bullet$ Julia Gaudio $^{a,\ast}$\\
        \small $^{a}$ Industrial Engineering and Management Sciences Department, Northwestern University, Evanston, IL, USA \\
        \small $^{b}$ Center for Engineering and Health, Northwestern University, Feinberg School of Medicine, Chicago, IL, USA \\\\
        \small $^{\ast}$Corresponding author: \tt{julia.gaudio@northwestern.edu} \\
}

\begin{document}

\maketitle

\normalsize
\begin{abstract}
In typical applications of facility location problems, the location of demand is assumed to be an input to the problem. The demand may be fixed or dynamic, but ultimately outside the optimizer's control. In contrast, there are settings, especially in humanitarian contexts, in which the optimizer decides where to locate a demand node. In this work, we introduce an optimization framework for joint facility and demand location. As examples of our general framework, we extend the well-known $k$-median and $k$-center problems into joint facility and demand location problems (JFDLP) and formulate them as integer programs. We propose a local search heuristic based on network flow. We apply our heuristic to a hurricane evacuation response case study. Our results demonstrate the challenging nature of these simultaneous optimization problems, especially when there are many potential locations. The local search heuristic is most promising when the the number of potential locations is large, while the number of facility and demand nodes to be located is small.
\end{abstract}

\section{Introduction}
The facility location problem is a classical problem in operations research, with numerous applications. There is an extensive literature on facility location models, algorithms, and applications \citep{Laporte2019,boonmee2017facility,Farahani2012,Melo2009,Farahani2009,Snyder2006,Drezner2004,Owen1998}. The facility location problem involves determining the optimal locations for facilities, where the precise objective depends on the application. For example, in supply chain networks, the locations of distribution centers are determined to minimize the cost of transporting goods to customers \citep{Amin2013,Melo2006,Daskin2005}. In urban planning, the locations of fire stations are chosen to maximize the number of people living within the serviceable area of the stations \citep{Drezner1991,Church1983}. In retail applications, a company may choose to open stores in order to maximize its market share \citep{Drezner2014,Drezner1998}. 

When formulating a facility location problem, the locations of the demand are typically considered to be intrinsic to the problem, and outside the optimizer's control. On the other hand, there are applications, particularly in the context of humanitarian and disaster relief, in which the locations of demand can be decided. For example, consider the problem of preparing for a natural disaster, which involves determining the locations of shelters for evacuees as well as the locations of warehouses to store aid supplies. Additionally, each shelter must be assigned to a warehouse, which may be capacitated. In this example, the warehouses would be considered facilities, and the shelters would be considered demand. In practice, these location decisions are typically made separately \citep{Sabbaghtorkan2020}, though there are several works proposing joint optimization of shelters and distribution centers to improve relief distribution \citep{Rodriguez2018,Dalal2018,Sheu2014}. For another example, consider refugee camp planning. While numerous factors influence the locations of tents in a camp, in arid regions, tent location decisions are strongly tied to water availability \citep{jahre2018approaches}. In these cases, the locations of wells and taps (facilities) must be jointly optimized with the locations of tents (demand) to ensure that each refugee has adequate access to water. For another example, consider the location of homeless shelters (demand) and facilities serving the homeless population, such as health clinics, family services, and housing agencies. Homeless people often travel great distances to meet their basic needs \citep{kaya2022improving,Wolch1993}; careful location of shelters and services has the potential to significantly improve service accessibility.

The common theme in the above examples is displacement, whether it be due to a natural disaster, a refugee crisis, or poverty. Displacement is becoming an increasingly common and urgent phenomenon, with more frequent and intense natural disasters \citep{landsea2005hurricanes}, recent geopolitical conflicts (i.e., Syrian Refugee Crisis (2011), Afghan Refugee Crisis (2021) and Ukraine-Russia War (2022)), and increasing rates of homelessness worldwide \citep{kuo2019yet}. Moreover, many of these joint location problems must be solved urgently, for example due to an incoming hurricane or an influx of refugees. With increasing displacement and changing conditions, humanitarian location problems are typically solved at the operational or tactical level of decision making \citep{balcik2008facility}. For example, refugee camps or shelter areas can be relocated according to the course of the disaster. In contrast, traditional facility location problems, such as factory or hospital location, correspond to strategic decisions, with more permanent consequences. With the increasing need for solving joint facility-demand location problems, we propose an optimization framework to capture these problems and design a heuristic to solve them.

In this paper, by formally conceptualizing \textit{demand location}, we introduce a new abstraction to the facility location literature, by allowing the optimizer to determine the location of the demand. To illustrate this new class of problems, consider the classical $k$-median problem. In this problem, there are $d$ fixed demand nodes. Given a set of $n$ possible facility locations, one must choose $k$ locations to open a facility. Each demand node is assigned to exactly one facility. The goal is to choose the locations of the $k$ facilities and assign the demand, in order to minimize the sum of assignment costs. Now suppose that the locations of the demand can also be chosen from among a set of $m$ possible locations. In this new problem, the optimizer must jointly decide the locations of the facilities, the locations of the demand, and the assignment of demand to facilities. In principle, any traditional facility location problem can be modified into a joint facility-demand location problem (JFDLP). Jointly optimizing the facility and demand locations is crucial in scenarios for which the dominating cost comes from transport between demand nodes and facility nodes. Finally, any JFDLP problem can include regional demand constraints, which impose lower- or upper-bounds on the number of demand nodes in each region.

We note that the addition of decision variables corresponding to the demand locations increases the complexity of the given facility location problem. Since facility location problems with fixed demand are themselves hard problems (many of them $\mathcal{NP}$-hard), we propose a local search heuristic. The heuristic is based on the following observation: when the facilities' locations are fixed, it is easy to determine the optimal demand locations. For example, consider the $k$-median problem with fixed facilities. In this case, the optimal demand locations can be determined by a network flow problem, which can be solved in polynomial time. This observation contrasts dramatically with the fact that the fixed-demand $k$-median problem is $\mathcal{NP}$-hard \citep{Hochbaum1984}. The key difference is the direction of the assignment: in the case of fixed facilities, each unit of demand is assigned to a single facility, while in the case of fixed demand, the assignments of facilities must cover the whole demand. 

In this paper, we illustrate the joint facility-demand location paradigm by extending two classical facility location problems: $k$-median and $k$-center. We introduce a local search algorithm for both problems, which iteratively updates the set of open facilities. Given a fixed set of facilities, we determine the optimal location of demand through a network flow problem, which can be solved in polynomial time. 

We evaluate the performance of the local search heuristic on a disaster management case study. The case study involves selecting hurricane shelters (demand) and warehouses (facilities) in Florida, under a $k$-median or $k$-center objective. The formulation includes regional constraints, prescribing lower bounds on the number of shelters in each region. We compare the performance of the local search algorithm to a commercial solver attempting to solve the exact problem. Our results show that the exact method perform poorly under the $k$-center objective and in some cases cannot even instantiate a solution. Though the heuristic is able to significantly improve upon the commercial solver when the number of potential locations is large, however the optimality is not guaranteed.  On the other hand, both methods perform well on the $k$-median objective, when the number of potential locations is small. The benefit of the heuristic over the commercial solver is seen in instances where the number of potential locations is large but the number of facilities and demands is small. Finally, we observe that the addition of regional constraints speeds up both the local search algorithm and the exact solver. Therefore, regional constraints may be added not just as a practical necessity, but also as a method to reduce the search space.

The remainder of this paper is organized as follows. The mathematical descriptions of our joint facility-demand location problems are presented in Section \ref{sec:formulation}. Section \ref{sec:algorithms} includes the local search heuristic, along with problem-dependent optimization-based subroutines. Section \ref{sec:case-study} describes the case study, with computational results appearing in Section \ref{sec:computational-results}. We discuss several directions for future work in Section \ref{sec:conclusion}.

\section{Mathematical Formulation}\label{sec:formulation}
We consider problems in which both the facilities and demand may be located among a set of discrete points. We consider the case of unit demand, where each demand point may accommodate exactly one unit of demand. Each demand unit must be assigned to a facility, where the facilities are capacitated. The goal is to optimize a given objective function, such as minimizing the sum of total assignment distances, or minimizing the largest assignment distance. Table \ref{table:sets} presents sets, decision variables and parameters used in integer programming formulation.

\begin{table}[H] \setlength{\tabcolsep}{4pt}
\caption{Sets, Decision Variables and Parameters}
\begin{center}
\begin{normalsize}
\begin{tabular}{cl}
\toprule
\textbf{Sets}&\textbf{Description}\\
\midrule
$\mathcal{I}$ & Set of candidate demand nodes, indexed by $i$, $|\mathcal{I}| = m$. \\
$\mathcal{J}$ & Set of candidate facility nodes, indexed by $j$, $|\mathcal{J}| = n$. \\
\toprule
\textbf{Decision Variables}&\textbf{Description}\\
\midrule
$x_{ij}$ & \textbf{1} if unit demand $i$ is assigned to facility $j$, \textbf{0} otherwise. \\
$y_{j}$ & \textbf{1} if a facility is located at candidate facility node $j$, \textbf{0} otherwise. \\
$z_{i}$ & \textbf{1} if unit demand is located at candidate demand node $i$, \textbf{0} otherwise. \\
\toprule
\textbf{Parameters}&\textbf{Description}\\
\midrule
$c_{ij}$ & Cost of traveling from candidate demand node $i$ to candidate facility node $j$. \\
$d$ & Total demand, satisfying $d \leq m$. \\
$k$ & Number of facilities to be opened, satisfying $k \leq n$. \\
$C$ & Capacity of facilities, where $Ck \geq d$. \\
\bottomrule
\end{tabular}
\end{normalsize}
\end{center}
\label{table:sets}
\end{table}

Here we present integer (binary) programs for the JFDLP-median and JFDLP-center problems.

\newpage

\textbf{Joint Facility and Demand Location Problem Under Median Objective:}

\begin{align}
    \text{min } & \sum_{i\in \mathcal{I}}\sum_{j\in \mathcal{J}}c_{ij}x_{ij} \label{eq_1} \\
    \text{s.t. } & \sum_{j\in \mathcal{J}}x_{ij} = z_i & \forall i \in \mathcal{I} \label{eq_2} \\
    & \sum_{j\in \mathcal{J}}y_{j} = k \label{eq_3} \\
    & \sum_{i\in \mathcal{I}}x_{ij} \leq Cy_{j} & \forall j \in \mathcal{J} \label{eq_4} \\
    & \sum_{i\in \mathcal{I}}z_{i} = d \label{eq_5} \\
    & x_{ij}, y_{j}, z_{i} \in \{0,1\} & \forall i \in \mathcal{I}, j \in \mathcal{J} \nonumber
\end{align}

In the joint facility and demand location problem under median objective, there are $n$ potential facility locations and $m$ potential demand locations. In this problem, we need to locate $k$ facilities and $d$ unit demands. Then, each demand node must be assigned to a facility, while minimizing the sum of distances. The objective function \ref{eq_1} minimizes the sum of total  distances from each unit demand node to the assigned facility.

Constraint \ref{eq_2} ensures that a unit demand node is assigned to only one facility. Constraint \ref{eq_3} makes sure that there are $k$ open facilities. Constraint \ref{eq_4} is the capacity constraint. Constraint \ref{eq_5} ensures that exactly $d$ units of demand are assigned. Note that if set of candidate demand nodes ($\mathcal{I}$) and set of candidate facility nodes ($\mathcal{J}$) are not disjoint, we can add the following constraint to prevent a demand and a facility from being located in the same place:

\begin{equation} \label{eq_31}
z_{j} + y_{j} \leq 1 \quad \quad \quad \forall j \in \mathcal{J}.
\end{equation}

The JFDLP under the center objective is expressed similarly.

\textbf{Joint Facility and Demand Location Problem Under Center Objective:}
\begin{align}
    \text{min } & \delta  \nonumber \\
    \text{s.t. } & c_{ij}x_{ij} \leq \delta & \forall i \in \mathcal{I}, j \in \mathcal{J} \label{eq_6} \\
    & \sum_{j\in \mathcal{J}}x_{ij} = z_i & \forall i \in \mathcal{I} \nonumber \\
    & \sum_{j\in \mathcal{J}}y_{j} = k \nonumber \\
    & \sum_{i\in \mathcal{I}}x_{ij} \leq Cy_{j} & \forall j \in \mathcal{J} \nonumber \\
    & \sum_{i\in \mathcal{I}}z_{i} = d \nonumber \\
    & x_{ij}, y_{j}, z_{i} \in \{0,1\} & \forall i \in \mathcal{I}, j \in \mathcal{J} \nonumber
\end{align}

In order to model the center objective, we introduce an auxiliary variable $\delta$. Due to Constraint \ref{eq_6}, any optimal solution will satisfy $\delta = \text{max}_{i \in \mathcal{I}} \{c_{ij}x_{ij}\}$.

\newpage

\subsection{Extension: Regional Constraints}

Any JFDLP problem can be modified to include regional constraints, which prescribe lower bounds, upper bounds, or equality constraints on the number of demand nodes in each (disjoint) region. Regional demand constraints allow the optimizer to incorporate real-world factors, which may either require a certain minimal level of demand, or may limit demand. These factors include road network infrastructure, disaster type (i.e., natural, human-made, slow-onset, sudden-onset), phase of disaster management cycle (i.e., mitigation, preparation, response, recovery), safety, security, aid distribution, access to healthcare, education and public services \citep{karsu2019refugee}. Additionally, regional constraints can be used as a computational tool, since they reduce the search space.

\begin{figure}
\begin{center}
\caption{Modeling of regional constraints as a flow problem.} \medskip
\includegraphics[scale=0.20]{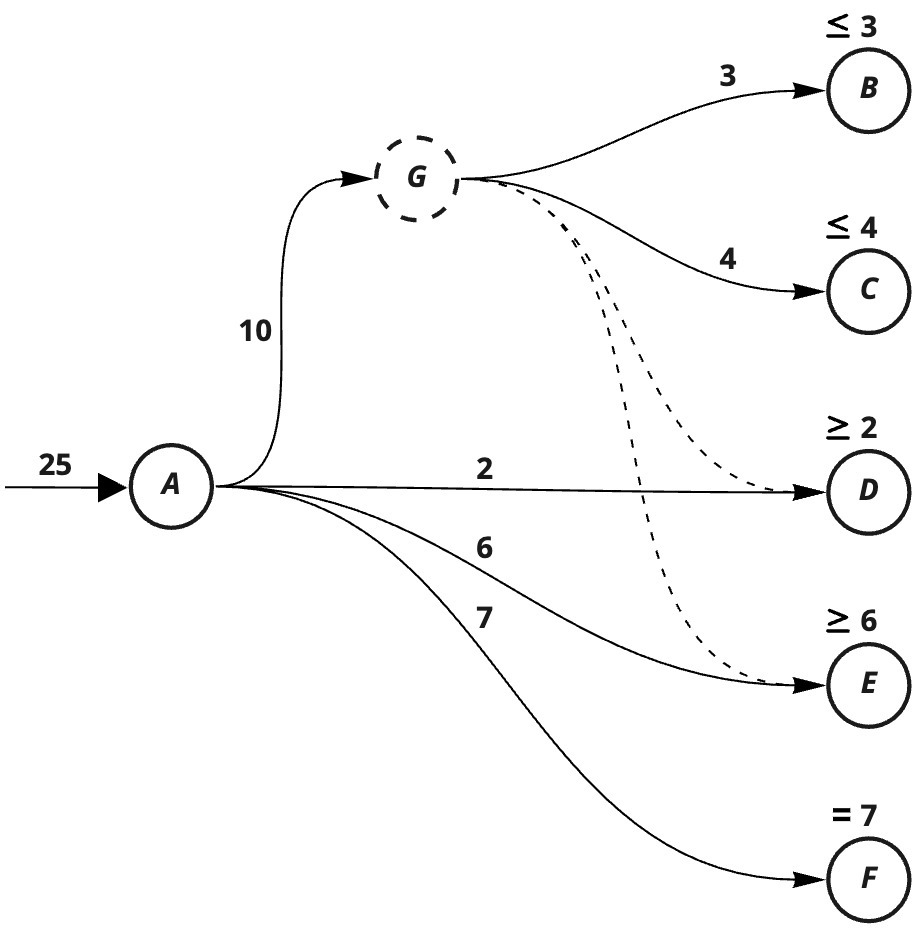}
\label{Img:Region-Flow}
\end{center}
\end{figure}

Regional constraints are modeled by standard network flow constructions. To illustrate the modeling of regional constraints, consider the case with five regions as in Figure \ref{Img:Region-Flow}. Node \textbf{A} is the source, sending 25 units of demand. Nodes \textbf{B}-\textbf{F} represent regional nodes, with constraints indicated above each node. These regional nodes are then connected downstream to the locations they represent (not pictured). We create an edge from \textbf{A} to the $\geq$ and $=$ constrained nodes, \textbf{D}, \textbf{E}, and \textbf{F}. The capacity of these edges matches the constraint value. The remaining flow $25 - (2+6+7) = 10$ is sent to a dummy node \textbf{G}. The dummy node is then connected to the $\leq$ constrained nodes, with capacities matching the regional capacities. Finally, the dummy node is connected to the $\geq$ constrained nodes. 

\FloatBarrier

\section{Algorithms}\label{sec:algorithms}

\begin{figure}
\begin{center}
\caption{Four types of demand/facility location problems with the $k$-median objective are illustrated. When both demand and facility nodes are fixed as in Figure (a), the optimal assignment of demand to facilities is a basic assignment problem, which can be expressed as a minimum cost network flow problem. When the demand is fixed but the facilities must be chosen (Figure (b)), we have the classic $k$-median facility location problem. On the other hand, when demand is variable but facilities are fixed (Figure (c)), we can optimize the $k$-median objective through a network flow problem. Finally, the JFDLP, illustrated in Figure (d), is not a network flow problem.} \medskip
\includegraphics[scale=0.40]{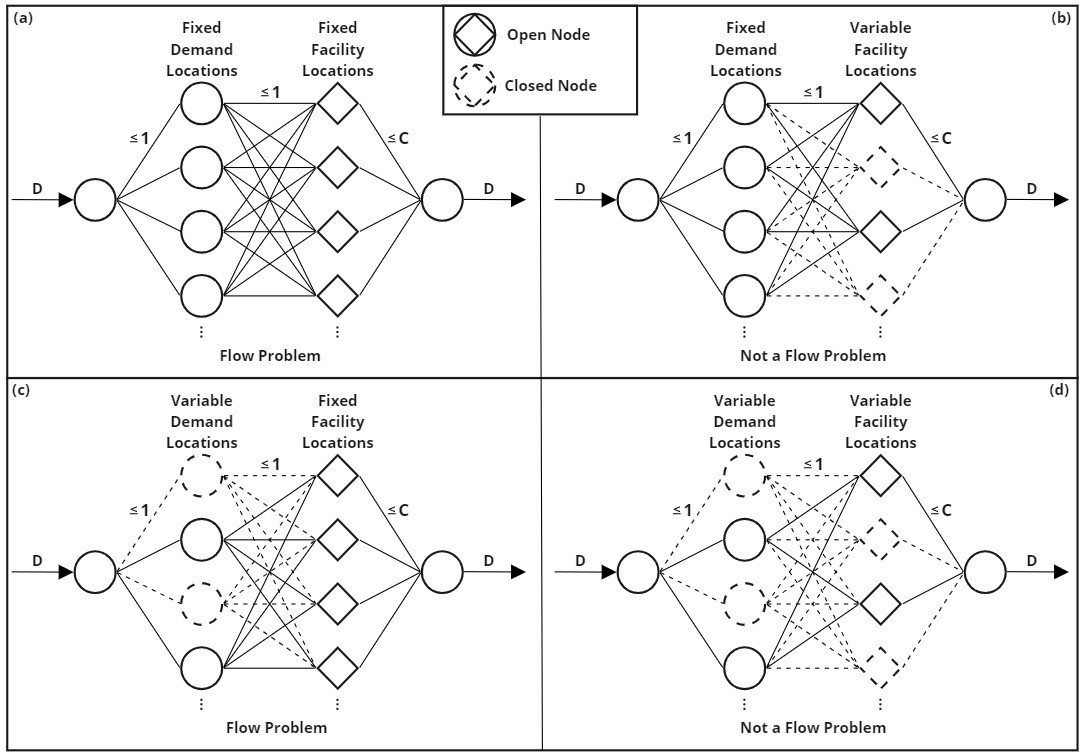}
\label{Image1}
\end{center}
\end{figure}

Our local search heuristic is based on the observation that when facilities are fixed, determining the optimal demand can be cast as a flow problem (Figure \ref{Image1}). The heuristic is described by the following algorithm. Let $\mathcal{F} = \{j|y_j = 1\}$ be the set of located facilities and let $\mathcal{D} = \{i|z_i = 1\}$ be the set of located demand nodes. 

\FloatBarrier

\begin{breakablealgorithm}
\caption{Local Search Algorithm}\label{alg:LSA}
\begin{algorithmic}[1]
\Require{Instance $I$ with parameters $m$, $n$, $d$, $k$, $C$, and costs $\{c_{ij}\}_{i \in \mathcal{I}, j \in \mathcal{J}}$}
\Ensure{Locations of facilities and demand}
\State $\mathcal{F} \gets \texttt{InitializeFacilities}(n,k)$
\State $obj \gets \texttt{FindOptimalDemand}(I, \mathcal{F})$

\State $(\text{OpenFacilities}, \text{ClosedFacilities}) \gets \texttt{RankFacilities}(\mathcal{F})$ \label{rank}
\For{$a \in \text{OpenFacilities}$}
\For{$b \in \text{ClosedFacilities}$}
\State $\overline{\mathcal{F}} \gets (\mathcal{F} \setminus\{a\}) \cup \{b\}$
\State $\overline{obj} \gets \texttt{FindOptimalDemand}(I, \overline{\mathcal{F}}, obj)$ 
\If{$\overline{obj} < obj$}
\State $obj \gets \overline{obj}$
\State $\mathcal{F} \gets \overline{\mathcal{F}}$
\State \Goto{rank}
\EndIf
\EndFor
\EndFor
\end{algorithmic}
\end{breakablealgorithm}
The initial choice of facility locations $\mathcal{F}$ in line $1$ is chosen uniformly at random from the set of size-$k$ subsets of the $n$ candidate facility locations. Given $\mathcal{F}$, the algorithm proposes a ``swap"; that is, opening a facility that is closed and closing a facility that is open. The new candidate solution, denoted $\overline{\mathcal{F}}$, is evaluated by $\texttt{FindOptimalDemand}$, which finds the optimal demand locations given the fixed facility locations $\overline{\mathcal{F}}$. The precise formulation of $\texttt{FindOptimalDemand}$ depends on the joint facility-demand location problem at hand. If the candidate solution improves on the current objective value, then we take another local search step, by returning to line $3$. The algorithm terminates when no local improvements can be made. We note that the order in which candidates are proposed influences how quickly an improvement is found. The $\texttt{RankFacilities}$ routine orders the open and closed facilities in a way that favors promising swaps. 

In order to solve JFDLP under the median objective, Algorithm \ref{alg:LSA} is used with Algorithm \ref{alg:Alg2}, where Algorithm \ref{alg:Alg2} is used in place of \texttt{FindOptimalDemand}. Since the optimization model is a network flow problem the integrality constraints are dropped and the model is formulated as a linear program.

\newpage

\begin{breakablealgorithm}
\caption{Fixed Facility Demand Location Problem With Median Objective}\label{alg:Alg2}
\begin{algorithmic}[1]
\Require{Instance $I$ with candidate sets $\mathcal{I}, \mathcal{J}$, parameters, $d$, $k$, $C$, costs $\{c_{ij}\}_{i \in \mathcal{I}, j \in \mathcal{J}}$ and open facilities $\mathcal{F}$} 
\Ensure{Optimal Objective Function Value, $x_{ij}$ and $z_{i}$}
\State Solve the following linear program

\begin{center}
\fbox{\begin{minipage}{25em}
\begin{align}
    \text{min } & \sum_{i\in \mathcal{I}}\sum_{j\in \mathcal{J}}c_{ij}x_{ij} \nonumber \\
    \text{s.t. } & \sum_{j\in \mathcal{J}}x_{ij} = z_i & \forall i \in \mathcal{I} \nonumber \\
    & \sum_{i\in \mathcal{I}}x_{ij} \leq C & \forall j \in \mathcal{F} \nonumber \\
    & \sum_{i\in \mathcal{I}}x_{ij} = 0 & \forall j \in \mathcal{J} \setminus \mathcal{F} \nonumber \\
    & \sum_{i\in \mathcal{I}}z_{i} = d \nonumber \\
    & 0 \leq x_{ij} \leq 1 & \forall i \in \mathcal{I}, j \in \mathcal{J} \nonumber \\
    & 0 \leq z_{i} \leq 1 & \forall i \in \mathcal{I} \nonumber
\end{align}
\end{minipage}}
\end{center}

\end{algorithmic}
\end{breakablealgorithm}

In order to solve the JFDLP under the center objective, Algorithm \ref{alg:Alg3} is used in place of \texttt{FindOptimalDemand}. This algorithm uses binary search with an LP. The approach is very similar to Algorithm \ref{alg:Alg2} but to accommodate the center objective, it employs a binary search, turning the optimization problem into a feasibility problem. Specifically, given a candidate distance $c_{\text{temp}}$, any edge with distance exceeding $c_{\text{temp}}$ has its capacity set to zero (Constraint \eqref{eq_25}).

\newpage

\begin{breakablealgorithm}
\caption{Fixed Facility Demand Location Problem With Center Objective}\label{alg:Alg3}
\begin{algorithmic}[1]
\Require{Instance $I$ with candidate sets $\mathcal{I}, \mathcal{J}$, parameters, $d$, $k$, $C$, costs $\{c_{ij}\}_{i \in \mathcal{I}, j \in \mathcal{J}}$ and open facilities $\mathcal{F}$, incumbent objective \textit{obj}}
\Ensure{Optimal $c_{max}$, $x_{ij}$ and $z_{i}$}
\State $\mathcal{S} \in \mathbb{R}^{m \times n} \gets$ values $\{c_{ij}: i \in \mathcal{I}, j \in \mathcal{J}\}$, sorted in increasing order.
\State $high \gets \text{Index of } obj \text{ in } \mathcal{S}$.
\State $low \gets 0$
\State $mid \gets \lfloor \frac{low + high}{2} \rfloor$
\While{$low < high - 1$}
\State $d_{temp} \gets \mathcal{S}$[$mid$]
\State Solve the following linear program

\begin{center}
\fbox{\begin{minipage}{25em}
\begin{align}
    \text{min } & \sum_{i\in \mathcal{I}}\sum_{j\in \mathcal{J}}c_{ij}x_{ij} \nonumber \\
    \text{s.t. } & \sum_{j\in \mathcal{J}}x_{ij} = z_i & \forall i \in \mathcal{I} \nonumber \\
    & \sum_{i\in \mathcal{I}}x_{ij} \leq C & \forall j \in \mathcal{F} \nonumber \\
    & \sum_{i\in \mathcal{I}}x_{ij} = 0 & \forall j \in \mathcal{J} \setminus \mathcal{F} \nonumber \\
    & \sum_{i\in \mathcal{I}}z_{i} = d \nonumber \\
    & x_{ij} = 0 & \forall i,j: c_{ij} > c_{temp} \label{eq_25} \\ 
    & 0 \leq x_{ij} \leq 1 & \forall i \in \mathcal{I}, j \in \mathcal{J} \nonumber \\
    & 0 \leq z_{i} \leq 1 & \forall i \in \mathcal{I} \nonumber
\end{align}
\end{minipage}}
\end{center}

\If{feasible}
\State $high \gets mid$
\Else
\State $low \gets mid + 1$
\EndIf
\EndWhile
\State Return $\mathcal{S}$[$mid$]
\end{algorithmic}
\end{breakablealgorithm}

\FloatBarrier

\begin{figure}[!hbt]
\begin{center}
\caption{Selection of an improving facility} \medskip
\includegraphics[scale=0.10]{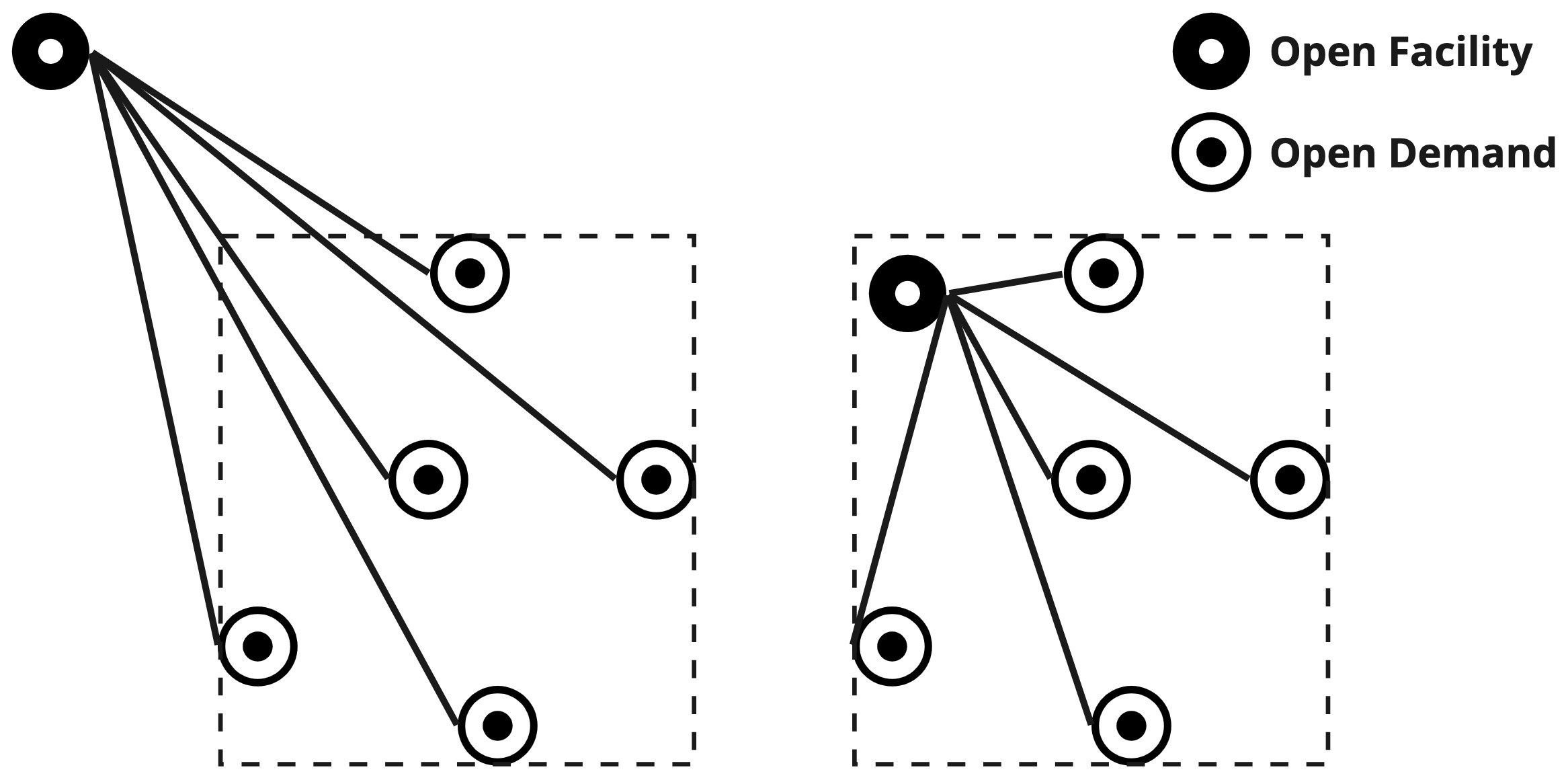}
\label{Img:ImprovementCases}
\end{center}
\end{figure}

Here we present an improvement heuristic algorithm (Algorithm \ref{alg:Alg4}) for the JFDLP under the center objective. Since the 1-center binary problem can be solved in polynomial time, we propose a polynomial time improvement heuristic algorithm, which splits a given solution into zones and solves the 1-center problem for each zone. Algorithm \ref{alg:Alg4} aims to eliminate solutions in the form depicted by Figure \ref{Img:ImprovementCases}, while maintaining feasibility.

The zone associated with facility location $j \in \mathcal{J}$ is defined as $\mathcal{Z}_j := \{i \in \mathcal{D} | x_{ij} = 1\}$. For each $j \in J$ such that $|\mathcal{Z}_j| > 0$, let $R_j$ be the smallest rectangle containing $\mathcal{Z}_j$. The improvement heuristic seeks to reassign all demand points in $\mathcal{Z}_j$ to a facility in $R_j$. By construction, all constraints remain satisfied. These reassignments are done sequentially across zones, while maintaining feasibility. We note that an improving facility for assigning the demand in zone $\mathcal{Z}_j$ may lie outside $R_j$. However, choosing such a facility may affect the improvement options of other zones; we therefore restrict the search for an improving facility to within $R_j$.

\FloatBarrier

\begin{breakablealgorithm}
\caption{Improvement Heuristic for JFDLP Under Center Objective}\label{alg:Alg4}
\begin{algorithmic}[1]
\Require{Locations of facilities $\mathcal{F}$, locations of demand $\mathcal{D}$, assignments $\{x_{ij}\}_{i \in \mathcal{I}, j \in \mathcal{J}}$}
\Ensure{Improved locations of facilities $\overline{\mathcal{F}}$}
\State Set $\overline{\mathcal{F}} = \emptyset$.
\State Determine $\mathcal{Z}_j = \{i \in \mathcal{D} \mid x_{ij} = 1\}$ for $j \in \mathcal{F}$.
\For{$j \in \mathcal{F}$} 
\State Determine the set of facilities locations within $R_j$, denoted $\mathcal{J}(R_j)$.
\State If $\mathcal{J}(R_j) \setminus \overline{\mathcal{F}} = \emptyset$, let $\overline{\mathcal{F}} = \mathcal{F}$ and \textbf{break}.
\State Find $j' \in \mathcal{J}(R_j) \setminus  \overline{\mathcal{F}}$ minimizing $\max_{i \in \mathcal{Z}_{j}} c_{ij'}$. 
\State If $\max_{i \in \mathcal{Z}_{j'}} c_{ij'} < \max_{i \in \mathcal{Z}_{j'}} c_{ij}$, set $\overline{F} = \overline{\mathcal{F}} \cup \{j'\}$. Otherwise, set $\overline{\mathcal{F}} = \overline{\mathcal{F}} \cup \{j\}$.
\EndFor
\end{algorithmic}
\end{breakablealgorithm}

\section{Case Study}\label{sec:case-study}

In this section, we describe our hurricane disaster management case study. The case study is set in Florida, which has been severely impacted by hurricanes in recent decades. Due to the sudden nature of hurricanes, decisions regarding evacuation, shelter, and resource prepositioning must be made quickly using the available information at hand. These decisions are influenced by hurricane predictions, infrastructure resilience, and resource availability. Public shelters play a critical role in disaster response. During the evacuation in response to Hurricane Irma in 2017, public shelters housed 190,000 evacuees \citep{Maul2018}. Still, a retrospective report found that more shelters were needed to meet the demand of evacuees, some of which did not evacuate because shelters were full \citep{Wong2018}. When we look at the past hurricane response practices, there were massive but few evacuation areas at long distances. This cause two problems. First, due to their risk perception evacuees try to reach to the evacuation areas using the same routes, causing extreme delays \citep{zhu2020estimating}. Second, the global COVID-19 pandemic has underscored the need for physical distancing in shelters \citep{nakai2021construction,pei2020compound,Shultz2020a}. At a given density level, larger shelters pose a greater health risk to evacuees than smaller shelters \citep{Shultz2020b}. As a response to these factors, our case study places many small hurricane shelters and supply warehouses.

The problem of locating shelters and the associated supply warehouses cleanly fits into the JFDLP framework. Here, the shelters are considered demand, while the warehouses are considered facilities. In this context, the median objective corresponds to minimizing the total transportation cost. On the other hand, the center objective minimizes the maximum distance between a shelter and its assigned supply warehouse, ensuring fairness \citep{huang2012models}.

\begin{figure}[ht!]
    \caption{GIS Data of City Centers (left) and ZIP Code Centers (right) for Florida}
    \label{fig:FLMaps}
    \includegraphics[width=.5\textwidth]{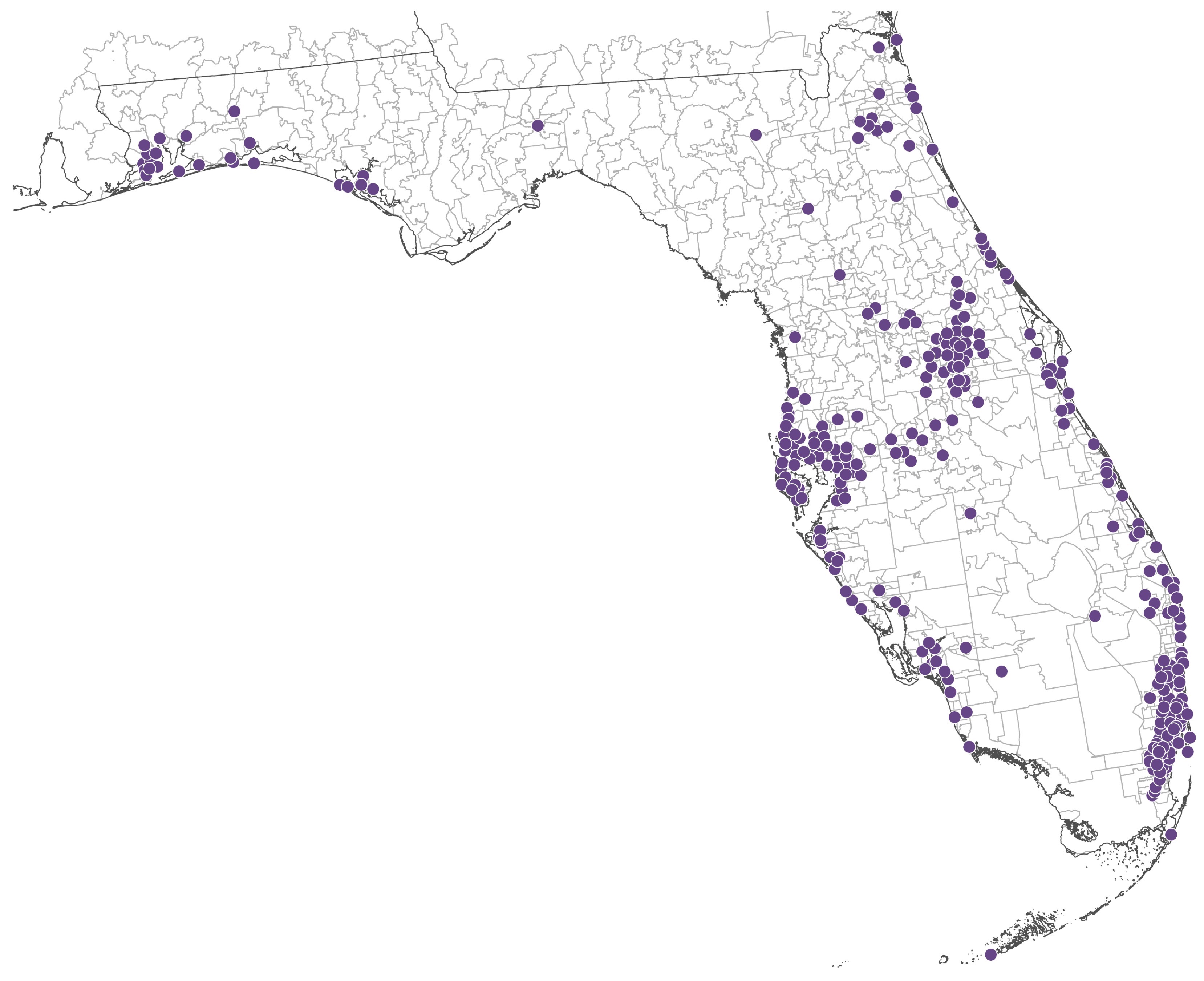}\hfill
    \includegraphics[width=.5\textwidth]{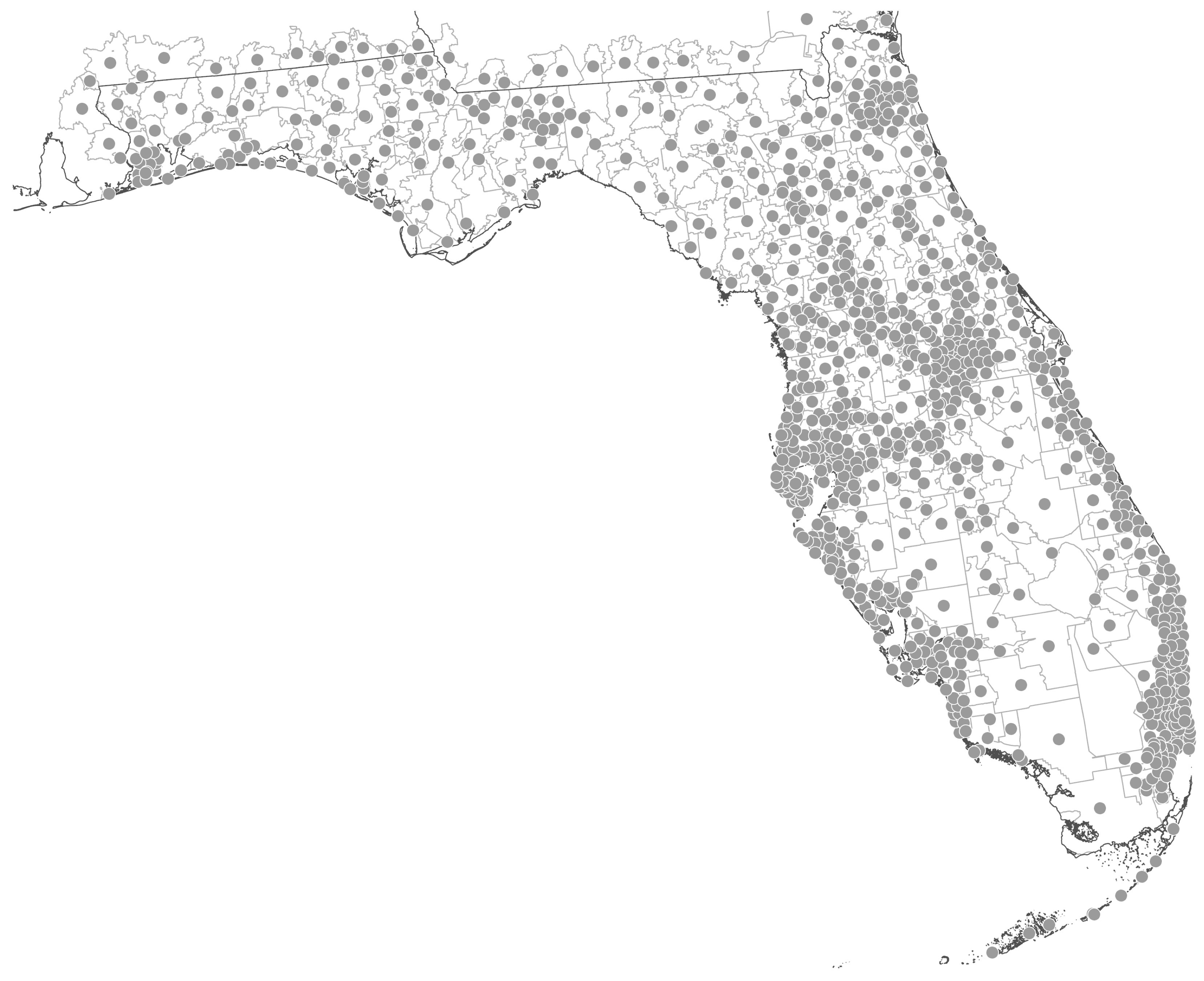}
\end{figure}

Figure \ref{fig:FLMaps} illustrates the city centers and ZIP code centers for Florida obtained using GIS. For each pair of locations $(i,j)$ the cost $c_{ij}$ is taken to be the travel time from $i$ to $j$, where $c_{ij}$ need not equal $c_{ji}$. The case study is done both for city centers and ZIP centers, in order to observe the effect of number of nodes on computational scaleability and complexity of exact methods and heuristic algorithms. There are 309 city centers and 1044 ZIP code centers in Florida. In City instances all of the cities are used, on the other hand for ZIP code instances, only the ZIP code centers that have a population more than 250 are included in the study, and this filtering reduces the number of ZIP code centers to 933. 

In our case study $\mathcal{I} = \mathcal{J}$. In order to prevent a demand and a facility node to be located at the same node we include Constraint \ref{eq_31} both in the exact and heuristic methods. 

\begin{table}[H] \setlength{\tabcolsep}{8pt}
\begin{center}
\caption{Values of Parameters in Different Sized Instances }
\begin{normalsize}
\begin{tabular}{lllll}
\toprule
\textbf{Notation} & \textbf{Description} & \multicolumn{3}{c}{\textbf{Instance Size}} \\
& & \textbf{Small}&\textbf{Medium}&\textbf{Large} \\
\midrule
\textit{d} & Number of binary demands to be located. & 50 & 100 & 150 \\
\textit{k} & Number of facilities to be located. & 3 & 6 & 9 \\
\textit{C} & Capacity of a facility. & 20 & 20 & 20 \\
\bottomrule
\end{tabular}
\end{normalsize}
\label{table:Instance-Size}
\end{center}
\end{table}

We consider multiple problem sizes, specified by Table \ref{table:Instance-Size}. We assume that each facility has the same capacity. In order to study the effect of regional partitions, we consider seven gridding configurations, after embedding the map of Florida into a square. Let $g$ be the grid parameter, where $g \in \{0,1,2,3,4,5,6\}$. For each $g$ value, the map of Florida is separated into $2^{2g}$ equal-sized non-overlapping squares. In order to add regional constraints, we define two new sets and a new set of parameters. Let $\mathcal{Q}$ be the set of non-overlapping regions, indexed by $q$ and let $\mathcal{R}_{q} = \mathcal{I} \cap q$ for $q \in \mathcal{Q}$. Finally, the parameter $\mu_q$ represents the required number of candidate demand nodes to be opened in region $q$.

In the context of hurricane preparedness and quick response, we only consider lower-bound regional constraint. In the case study to approximate $\mu_q$ parameter, 70\% of the overall demand $d$ is distributed across the grid, proportional to the grid population (rounding down), which was inspired by existing practices in \citet{wilmot2005methodology}. The following constraint is included in both the exact and heuristic models:

\begin{equation*}
\sum_{i\in R_q}z_{i} \geq \mu_q \quad \quad \quad \forall q \in Q.
\end{equation*}

\section{Computational Results}\label{sec:computational-results}

In this section we present the computational results for our experiments. We compare the objective values of the exact and heuristic methods, and report the percent difference:

\begin{equation} \nonumber
\text{Percent Difference} = \frac{\text{Heuristic Objective Function Value - Exact Objective Function Value}}{\text{Exact Objective Function Value}} \times 100\%
\end{equation}

All of the models are solved by Gurobi 9.0.3 on a desktop computer with an Intel(R) Core (TM) i7-9750H CPU @ 2.60GHz processor and a 16GB RAM.

Tables \ref{table:results-city-median}-\ref{table:results-ZIP-center} present the experimental results for both city and ZIP code instances, under both the median and center objectives and Table \ref{table:improvement} presents the results for ZIP codes under the center objective using the improvement heuristic (Algorithm \ref{alg:Alg4}). In city instances $m=n=309$ and in ZIP instances $m=n=933$. In the tables, ``Counter'' refers to the number of times that the local search algorithm (Algorithm \ref{alg:LSA}) executes \textit{line} 11. The MIP gap is the percentage difference between MIP objective bound and incumbent solution objective. 

\begin{equation} \nonumber
\text{MIP Gap} = \frac{|\text{Incumbent Objective Value - Dual Objective Bound}|}{\text{Incumbent Objective Value}} \times 100\%
\end{equation}

Finally if a run time cell has asterisk sign (*) it means that the method has been stopped due to exceeding the time limit of one hour.

\begin{table}[!hbt] \setlength{\tabcolsep}{8pt}
\begin{center}
\caption{Computational Results of JFDLP Under Median Objective for Cities}
\begin{adjustbox}{width=0.9\textwidth}
\begin{tabular}{ll|lll|lll|l}
\toprule
\multicolumn{2}{c}{\textbf{Instance}} & \multicolumn{3}{c}{\textbf{Exact Method}} & \multicolumn{3}{c}{\textbf{Heuristic Method}} & \textbf{Percentage}\\
\textbf{Size} & \textbf{Grid} & \textbf{Run Time (s)} & \textbf{Objective Value} & \textbf{MIP Gap} & \textbf{Run Time (s)} & \textbf{Objective Value} & \textbf{Counter} & \textbf{Difference}\\
\midrule
S & 0 & 20.15 & 606.70 & 0.0\% & 82.85 & 606.70 & 9 & 0.00\% \\
S & 1 & 50.79 & 699.46 & 0.0\% & 128.56 & 699.46 & 10 & 0.00\% \\
S & 2 & 2325.59 & 1209.65 & 0.0\% & 166.20 & 1209.65 & 15 & 0.00\% \\
S & 3 & 2358.28 & 980.30 & 0.0\% & 180.77 & 980.30 & 17 & 0.00\% \\
S & 4 & 1301.57 & 996.89 & 0.0\% & 185.41 & 996.89 & 13 & 0.00\% \\
S & 5 & 143.89 & 841.71 & 0.0\% & 150.68 & 841.71 & 10 & 0.00\% \\
\hline
M & 0 & 22.90 & 1373.58 & 0.0\% & 867.43 & 1373.58 & 16 & 0.00\% \\
M & 1 & 55.89 & 1405.37 & 0.0\% & 1253.40 & 1405.37 & 20 & 0.00\% \\
M & 2 & 230.34 & 1773.10 & 0.0\% & 1933.28 & 1773.10 & 40 & 0.00\% \\
M & 3 & 1459.61 & 1961.32 & 0.0\% & 1766.04 & 1961.32 & 34 & 0.00\% \\
M & 4 & 1009.76 & 1995.48 & 0.0\% & 1794.48 & 1995.48 & 32 & 0.00\% \\
M & 5 & 176.21 & 1675.19 & 0.0\% & 613.03 & 1675.19 & 18 & 0.00\% \\
\hline
L & 0 & 93.95 & 2308.41 & 0.0\% & 2276.43 & 2308.41 & 24 & 0.00\% \\
L & 1 & 76.70 & 2308.41 & 0.0\% & 3600.11* & 2372.57 & 33 & 2.78\% \\
L & 2 & 128.52 & 2591.37 & 0.0\% & 2591.37 & 2591.37 & 33 & 0.00\% \\
L & 3 & 2750.53 & 2929.09 & 0.0\% & 3600.07* & 3049.69 & 37 & 4.12\% \\
L & 4 & 1098.90 & 2959.68 & 0.0\% & 3600.10* & 3506.02 & 39 & 18.46\% \\
L & 5 & 184.51 & 2704.47 & 0.0\% & 3007.13 & 2704.47 & 28 & 0.00\% \\
\bottomrule
\end{tabular}
\end{adjustbox}
\label{table:results-city-median}
\end{center}
\end{table}

\begin{table}[!hbt] \setlength{\tabcolsep}{8pt}
\begin{center}
\caption{Computational Results of JFDLP Under Median Objective for ZIPs}
\begin{adjustbox}{width=0.9\textwidth}
\begin{tabular}{ll|lll|lll|l}
\toprule
\multicolumn{2}{c}{\textbf{Instance}} & \multicolumn{3}{c}{\textbf{Exact Method}} & \multicolumn{3}{c}{\textbf{Heuristic Method}} & \textbf{Percentage}\\
\textbf{Size} & \textbf{Grid} & \textbf{Run Time (s)} & \textbf{Objective Value} & \textbf{MIP Gap} & \textbf{Run Time (s)} & \textbf{Objective Value} & \textbf{Counter} & \textbf{Difference}\\
\midrule
S & 0 & 382.04 & 518.90 & 0.0\% & 1270.21 & 518.90 & 8 & 0.00\% \\
S & 1 & 506.56 & 570.76 & 0.0\% & 1025.01 & 570.76 & 17 & 0.00\% \\
S & 2 & 3636.37* & 1123.34 & 48.8\% & 1123.34 & 1123.34 & 13 & 0.00\% \\
S & 3 & 3636.70* & 921.95 & 38.9\% & 1997.45 & 918.22 & 18 & -0.40\% \\
S & 4 & 3636.63* & 774.17 & 31.1\% & 1856.94 & 741.21 & 16 & -4.26\% \\
S & 5 & 461.38 & 555.13 & 0.0\% & 1278.61 & 555.13 & 14 & 0.00\% \\
\hline
M & 0 & 270.39 & 1086.14 & 0.0\% & 3600.21* & 1132.42 & 16 & 4.26\% \\
M & 1 & 332.75 & 1124.00 & 0.0\% & 3600.09* & 1218.55 & 16 & 8.41\% \\
M & 2 & 1081.50 & 1244.09 & 0.0\% & 3600.18* & 1435.92 & 14 & 15.42\% \\
M & 3 & 3637.17* & 1563.36 & 23.6\% & 3600.13* & 1845.27 & 20 & 18.03\% \\
M & 4 & 3637.39* & 1376.98 & 13.2\% & 3600.08* & 1515.25 & 16 & 10.04\% \\
M & 5 & 3636.67* & 1192.44 & 3.0\% & 3600.15* & 1447.90 & 16 & 21.42\% \\
\hline
L & 0 & 440.25 & 1681.19 & 0.0\% & 3600.01* & 2039.84 & 16 & 21.33\% \\
L & 1 & 379.67 & 1718.75 & 0.0\% & 3600.13* & 1952.80 & 13 & 13.62\% \\
L & 2 & 1356.46 & 1849.68 & 0.0\% & 3600.29* & 2605.08 & 9 & 40.84\% \\
L & 3 & 3636.44* & 2174.83 & 13.1\% & 3600.31* & 2905.69 & 16 & 33.61\% \\
L & 4 & 3637.30* & 2252.32 & 16.9\% & 3600.01* & 2777.69 & 8 & 23.33\% \\
L & 5 & 3638.05* & 1822.97 & 2.9\% & 3600.23* & 2358.72 & 19 & 29.39\% \\
\bottomrule
\end{tabular}
\end{adjustbox}
\label{table:results-ZIP-median}
\end{center}
\end{table}

\begin{table}[!hbt] \setlength{\tabcolsep}{8pt}
\begin{center}
\caption{Computational Results of JFDLP Under Center Objective for Cities}
\begin{adjustbox}{width=0.9\textwidth}
\begin{tabular}{ll|lll|lll|l}
\toprule
\multicolumn{2}{c}{\textbf{Instance}} & \multicolumn{3}{c}{\textbf{Exact Method}} & \multicolumn{3}{c}{\textbf{Heuristic Method}} & \textbf{Percentage}\\
\textbf{Size} & \textbf{Grid} & \textbf{Run Time (s)} & \textbf{Objective Value} & \textbf{MIP Gap} & \textbf{Run Time (s)} & \textbf{Objective Value} & \textbf{Counter} & \textbf{Difference}\\
\midrule
S & 0 & 1111.33 & 17.80 & 0.0\% & 1336.94 & 22.80 & 10 & 28.06\% \\
S & 1 & 1123.24 & 26.64 & 0.0\% & 3600.02* & 50.94 & 22 & 91.12\% \\
S & 2 & 1460.72 & 96.51 & 0.0\% & 3600.31* & 152.05 & 23 & 57.55\% \\
S & 3 & 3605.59* & 56.68 & 24.2\% & 1549.27 & 101.42 & 7 & 79.93\% \\
S & 4 & 3605.17* & 54.06 & 1.4\% & 3600.34* & 76.22 & 15 & 40.99\% \\
S & 5 & 3605.12* & 63.00 & 37.2\% & 3600.38* & 85.85 & 14 & 36.26\% \\
\hline
M & 0 & 892.33 & 21.41 & 0.0\% & 3600.28* & 32.99 & 5 & 54.07\% \\
M & 1 & 1467.21 & 21.92 & 0.0\% & 3600.53* & 109.36 & 12 & 398.90\% \\
M & 2 & 158.53 & 60.41 & 0.0\% & 3600.95* & 110.14 & 7 & 82.32\% \\
M & 3 & 3609.97* & 60.24 & 26.7\% & 3600.77* & 191.16 & 17 & 217.33\% \\
M & 4 & 3605.47* & 64.60 & 35.4\% & 3601.26* & 168.23 & 11 & 160.43\% \\
M & 5 & 78.94 & 38.29 & 0.0\% & 3600.19* & 90.72 & 8 & 137.91\% \\
\hline
L & 0 & 2285.61 & 25.62 & 0.0\% & 3600.58* & 38.48 & 5 & 50.22\% \\
L & 1 & 2135.59 & 26.24 & 0.0\% & 3600.72* & 148.48 & 12 & 465.85\% \\
L & 2 & 183.53 & 54.42 & 0.0\% & 3600.40* & 208.98 & 15 & 284.02\% \\
L & 3 & 13.74 & 102.59 & 0.0\% & 3600.48* & 131.22 & 6 & 27.91\% \\
L & 4 & 29.70 & 102.59 & 0.0\% & 3600.63* & 135.96 & 6 & 32.53\% \\
L & 5 & 15.00 & 102.59 & 0.0\% & 3600.45* & 102.59 & 4 & 0.00\% \\
\bottomrule
\end{tabular}
\end{adjustbox}
\label{table:results-city-center}
\end{center}
\end{table}

\begin{table}[!hbt] \setlength{\tabcolsep}{8pt}
\begin{center}
\caption{Computational Results of JFDLP Under Center Objective for ZIPs}
\begin{adjustbox}{width=0.9\textwidth}
\begin{tabular}{ll|lll|lll|l}
\toprule
\multicolumn{2}{c}{\textbf{Instance}} & \multicolumn{3}{c}{\textbf{Exact Method}} & \multicolumn{3}{c}{\textbf{Heuristic Method}} & \textbf{Percentage}\\
\textbf{Size} & \textbf{Grid} & \textbf{Run Time (s)} & \textbf{Objective Value} & \textbf{MIP Gap} & \textbf{Run Time (s)} & \textbf{Objective Value} & \textbf{Counter} & \textbf{Difference}\\
\midrule
S & 0 & 3646.92* & 449.12 & 100.0\% & 3601.18* & 32.86 & 1 & -92.68\% \\
S & 1 & 3646.93* & 449.12 & 100.0\% & 3601.68* & 128.02 & 13 & -71.50\% \\
S & 2 & 3646.15* & 457.80 & 100.0\% & 3600.97* & 227.31 & 13 & -50.35\% \\
S & 3 & 3646.27* & 421.00 & 100.0\% & 3601.31* & 111.07 & 3 & -73.62\% \\
S & 4 & 3646.83* & 378.96 & 100.0\% & 3601.90* & 142.36 & 3 & -62.43\% \\
S & 5 & 3646.49* & 368.61 & 100.0\% & 3602.34* & 78.41 & 2 & -78.73\% \\
\hline
M & 0 & 3646.52* & 449.12 & 100.0\% & 3603.33* & 32.60 & 1 & -92.74\% \\
M & 1 & 3646.68* & 609.72 & 100.0\% & 3603.23* & 59.73 & 2 & -90.20\% \\
M & 2 & 3646.19* & 480.01 & 100.0\% & 3600.94* & 109.38 & 1 & -77.21\% \\
M & 3 & 3646.44* & 457.49 & 100.0\% & 3602.12* & 209.49 & 2 & -54.21\% \\
M & 4 & 3646.85* & 457.80 & 100.0\% & 3602.37* & 167.95 & 3 & -63.31\% \\
M & 5 & 3646.54* & 185.83 & 100.0\% & 3600.18* & 159.80 & 1 & -14.01\% \\
\hline
L & 0 & 3646.26* & 606.42 & 100.0\% & 3605.73* & 23.91 & 1 & -96.06\% \\
L & 1 & 3646.86* & 640.60 & 100.0\% & 3602.59* & 136.38 & 2 & -78.71\% \\
L & 2 & 3646.26* & 574.42 & 100.0\% & 3606.08* & 80.52 & 1 & -85.98\% \\
L & 3 & 3646.65* & 563.44 & 100.0\% & 3601.76* & 119.07 & 1 & -78.87\% \\
L & 4 & 3646.81* & 606.42 & 100.0\% & 3603.20* & 121.06 & 1 & -80.04\% \\
L & 5 & 3646.21* & 560.33 & 100.0\% & 3606.39* & 97.97 & 2 & -82.52\% \\
\bottomrule
\end{tabular}
\end{adjustbox}
\label{table:results-ZIP-center}
\end{center}
\end{table}

\FloatBarrier

Examining the results of Table \ref{table:results-city-median} (median objective, cities), the exact algorithm solves the problem to optimality in all of the instances, and the heuristic achieves the optimal solution in all but three instances. The heuristic is significantly faster for some small instances, while the exact method dominates larger instances. When we consider the median objective with ZIP code locations (Table \ref{table:results-ZIP-median}), the exact method fails to solve the problem to optimality for many of the instances, particularly when $g \geq 2$. The heuristic outperforms the exact method in smaller instances, but is generally outperformed for larger instances. 

When we look at the solutions of JFDLP under the center objective applied to cities (Table \ref{table:results-city-center}), the exact method generally solves all instances to optimality. In five cases, it hits the time limit and reports a non-zero MIP gap. On the other hand, the heuristic method reports large percentage differences in nearly all instances. Finally, Table \ref{table:results-ZIP-center} presents the results for the center objective with ZIP code locations. Both the exact and heuristic methods hit the one hour time limit in all of the cases. The exact solver reports a 100\% MIP gap in all instances, while the heuristic is able to improve significantly relative to the exact method. 

Since ZIP codes are spread uniformly (Figure \ref{fig:FLMaps}), we apply the improvement heuristic in this setting. Table \ref{table:improvement} summarizes the results. The algorithm runs less than one second, giving an average improvement of -23.41\%. Only in non-gridded instances does the improvement method not result in an improvement, due to overlapping zones. 

\begin{equation} \nonumber
\text{Percent difference} = \frac{\text{New Objective Function Value - Old Objective Function Value}}{\text{Old Objective Function Value}} \times 100\%
\end{equation}

\FloatBarrier

\begin{table}[!hbt] \setlength{\tabcolsep}{8pt}
\begin{center}
\caption{Computational Results of JFDLP Under Center Objective for ZIP Codes After Improvement Heuristic}
\begin{adjustbox}{width=0.6\textwidth}
\begin{tabular}{ll|l|ll|l}
\toprule
\multicolumn{2}{c}{\textbf{Instance}} & \multicolumn{1}{c}{\textbf{Heuristic Method}} & \multicolumn{2}{c}{\textbf{Improvement Heuristic Method}} & \textbf{Percentage}\\
\textbf{Size} & \textbf{Grid} & \textbf{Objective Value} & \textbf{Objective Value} & \textbf{Run Time (s)}  & \textbf{Difference}\\
\midrule
S & 0 & 32.86 & 32.86 & 0.03 & 00.00\% \\
S & 1 & 128.02 & 84.48 & 0.16 & -34.01\% \\
S & 2 & 227.31 & 147.71 & 0.29 & -35.02\% \\
S & 3 & 111.07 & 100.80 & 0.10 & -9.25\% \\
S & 4 & 142.36 & 96.74 & 0.09 & -32.05\% \\
S & 5 & 78.41 & 52.78 & 0.05 & -32.69\% \\
\hline
M & 0 & 32.60 & 35.15 & 0.07 & 0.00\% \\
M & 1 & 59.73 & 41.31 & 0.05 & -30.84\% \\
M & 2 & 109.38 & 86.11 & 0.10 & -21.27\% \\
M & 3 & 209.49 & 119.17 & 0.24 & -43.11\% \\
M & 4 & 167.95 & 108.02 & 0.18 & -35.68\% \\
M & 5 & 159.80 & 116.00 & 0.11 & -27.41\% \\
\hline
L & 0 & 23.91 & 25.88 & 0.05 & 0.00\% \\
L & 1 & 136.38 & 90.36 & 0.08 & -33.74\% \\
L & 2 & 80.52 & 62.12 & 0.04 & -22.85\% \\
L & 3 & 119.07 & 109.02 & 0.12 & -8.44\% \\
L & 4 & 121.06 & 88.05 & 0.19 & -27.27\% \\
L & 5 & 97.97 & 70.87 & 0.15 & -27.66\% \\
\bottomrule
\end{tabular}
\end{adjustbox}
\label{table:improvement}
\end{center}
\end{table}

\FloatBarrier

\begin{figure}[h]
    \caption{Performance Comparison of Heuristic and Exact Method for the Small Instance and Grid Size 2}
    \label{Fig:PerformanceComparison}
     \centering
     \begin{subfigure}[h]{0.45\textwidth} 
         \centering
         \caption{Instance = City, Objective = Median} \label{SubFig:a}
         \includegraphics[width=0.80\textwidth]{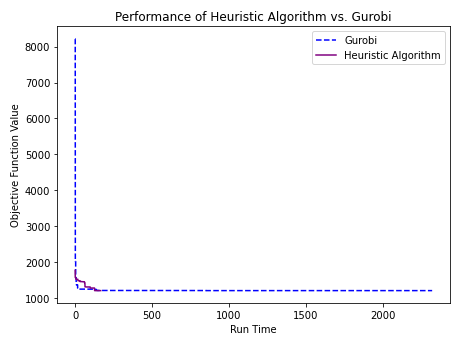}
     \end{subfigure}
     \hfill
     \begin{subfigure}[h]{0.45\textwidth} 
         \centering
         \caption{Instance = ZIP, Objective = Median} \label{SubFig:b}
         \includegraphics[width=0.80\textwidth]{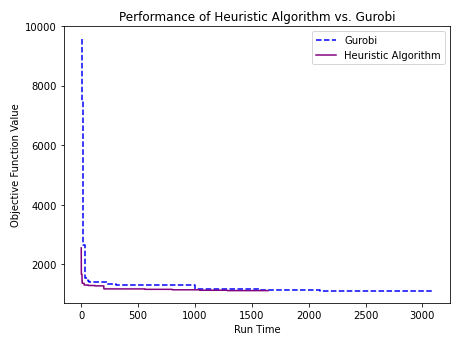}
     \end{subfigure} \\
     \centering
     \begin{subfigure}[h]{0.45\textwidth} 
         \centering
         \caption{Instance = City, Objective = Center} \label{SubFig:c}
         \includegraphics[width=0.80\textwidth]{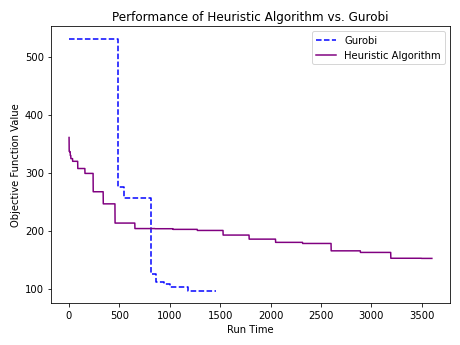}
     \end{subfigure}
     \hfill
     \begin{subfigure}[h]{0.45\textwidth}
         \centering
         \caption{Instance = ZIP, Objective = Center} \label{SubFig:d}
         \includegraphics[width=0.80\textwidth]{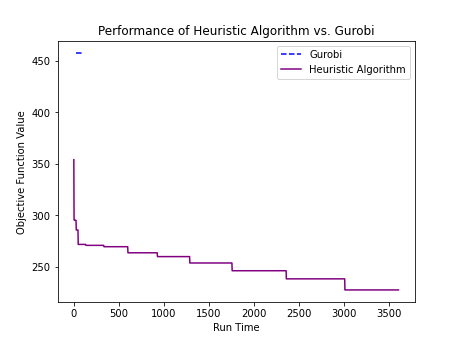}
     \end{subfigure}
\end{figure}

\FloatBarrier

Figure \ref{Fig:PerformanceComparison} compares the performance trajectories of the exact and heuristic methods for small instance and grid size 2, corresponding to third line of Tables \ref{table:results-city-median} to \ref{table:results-ZIP-center}. For the median objective, the heuristic algorithm is much faster than the exact method both in city and ZIP codes (Figures \ref{SubFig:a} and \ref{SubFig:b}). Under the center objective (cities), the heuristic algorithm performs better at the beginning, but then the exact algorithm passes the heuristic at approximately $t =700s$ (Figure \ref{SubFig:c}). When locations are given by ZIP codes, the commercial solver cannot start solving the problem but our heuristic method continuously finds improving solutions (Figure \ref{SubFig:d}).

\FloatBarrier

\begin{figure}[h]
    \caption{Visualization of Solutions to the Small Instance and Grid Size 2 Under Median Objective}
    \label{Fig:MedianMap}
     \centering
     \begin{subfigure}[h]{0.45\textwidth}
         \centering
         \caption{Instance = City, Method = Exact}
         \includegraphics[width=0.60\textwidth]{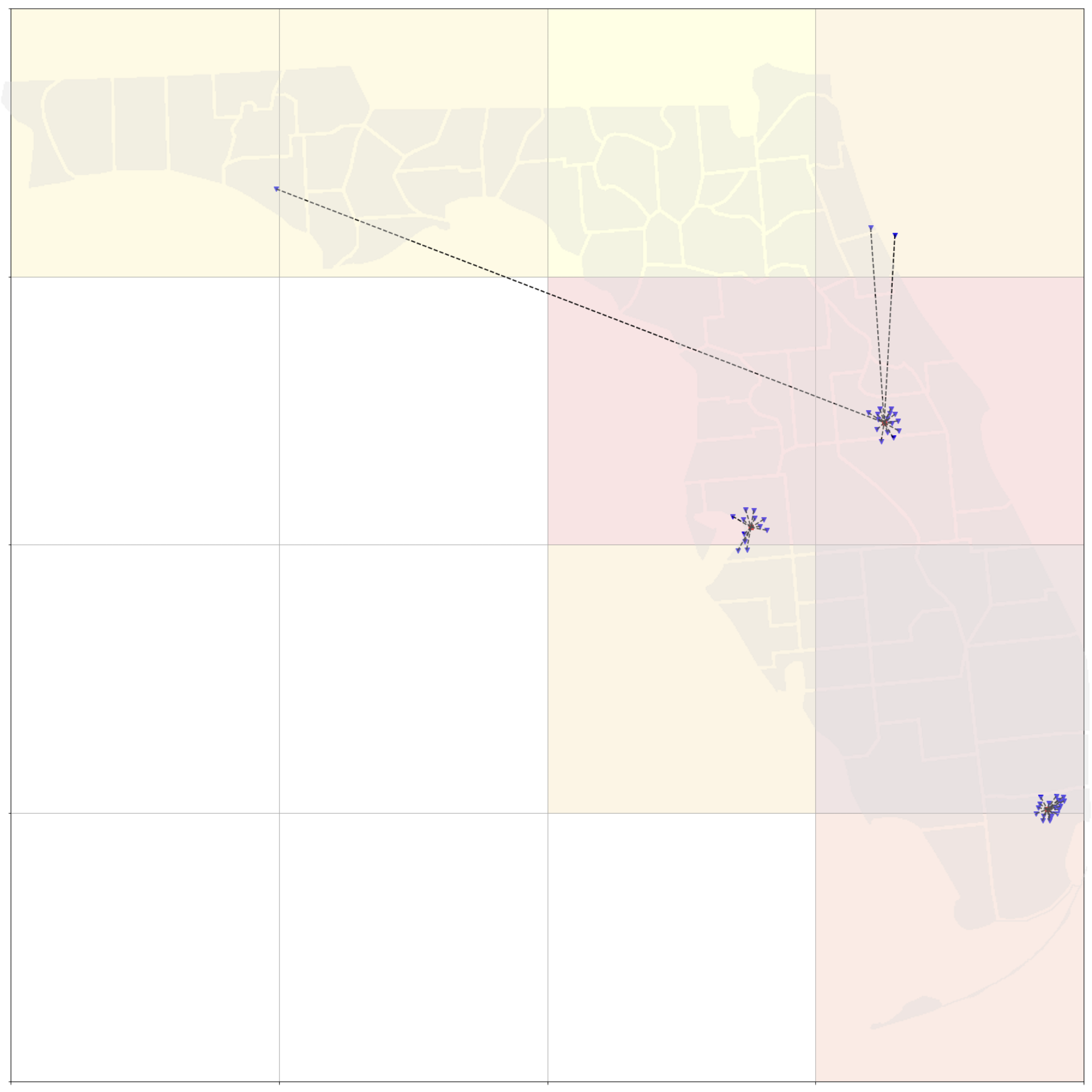}
     \end{subfigure}
     \hfill
     \begin{subfigure}[h]{0.45\textwidth}
         \centering
         \caption{Instance = City, Method = Heuristic}
         \includegraphics[width=0.60\textwidth]{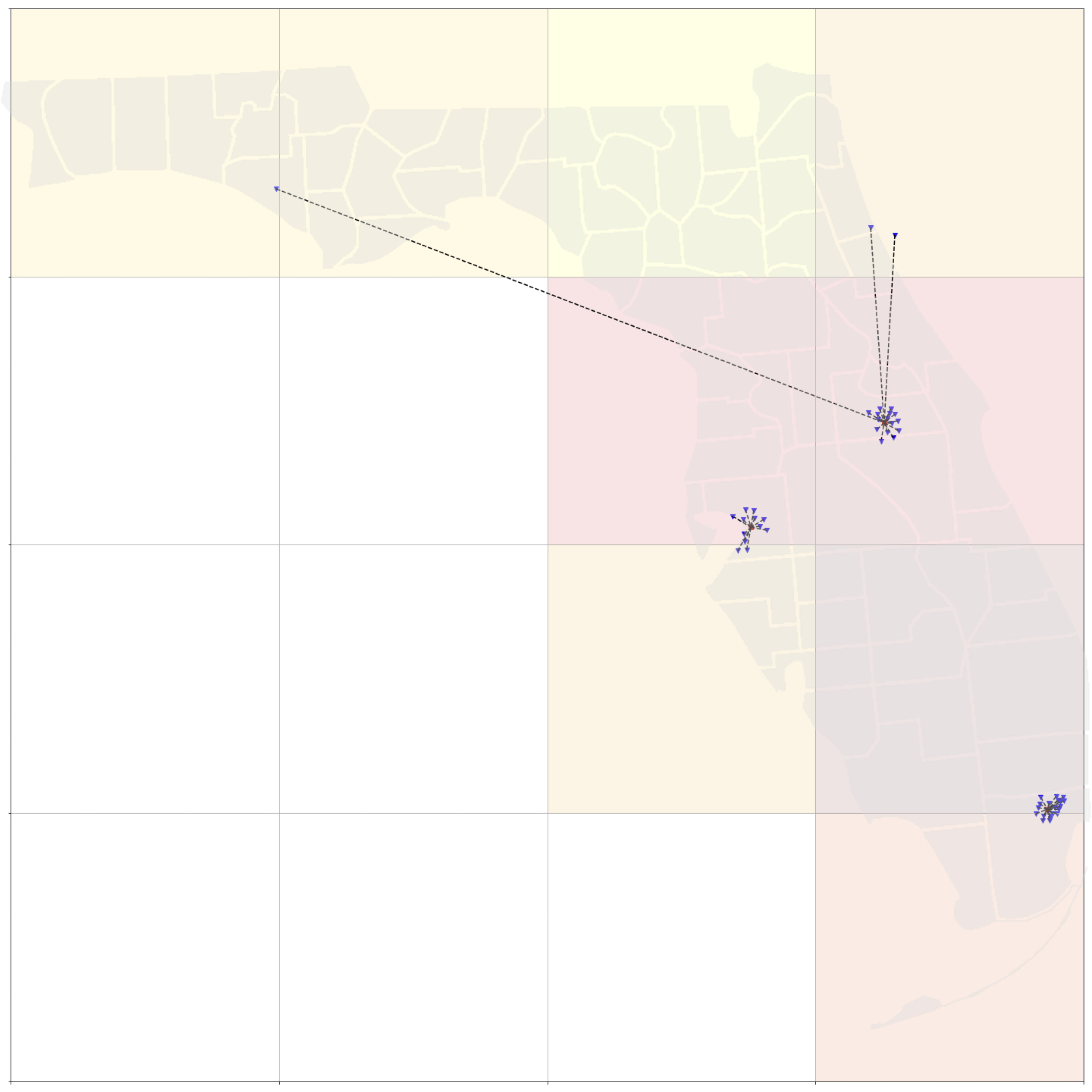}
     \end{subfigure} \\
     \centering
     \begin{subfigure}[h]{0.45\textwidth}
         \centering
         \caption{Instance = ZIP, Method = Exact}
         \includegraphics[width=0.60\textwidth]{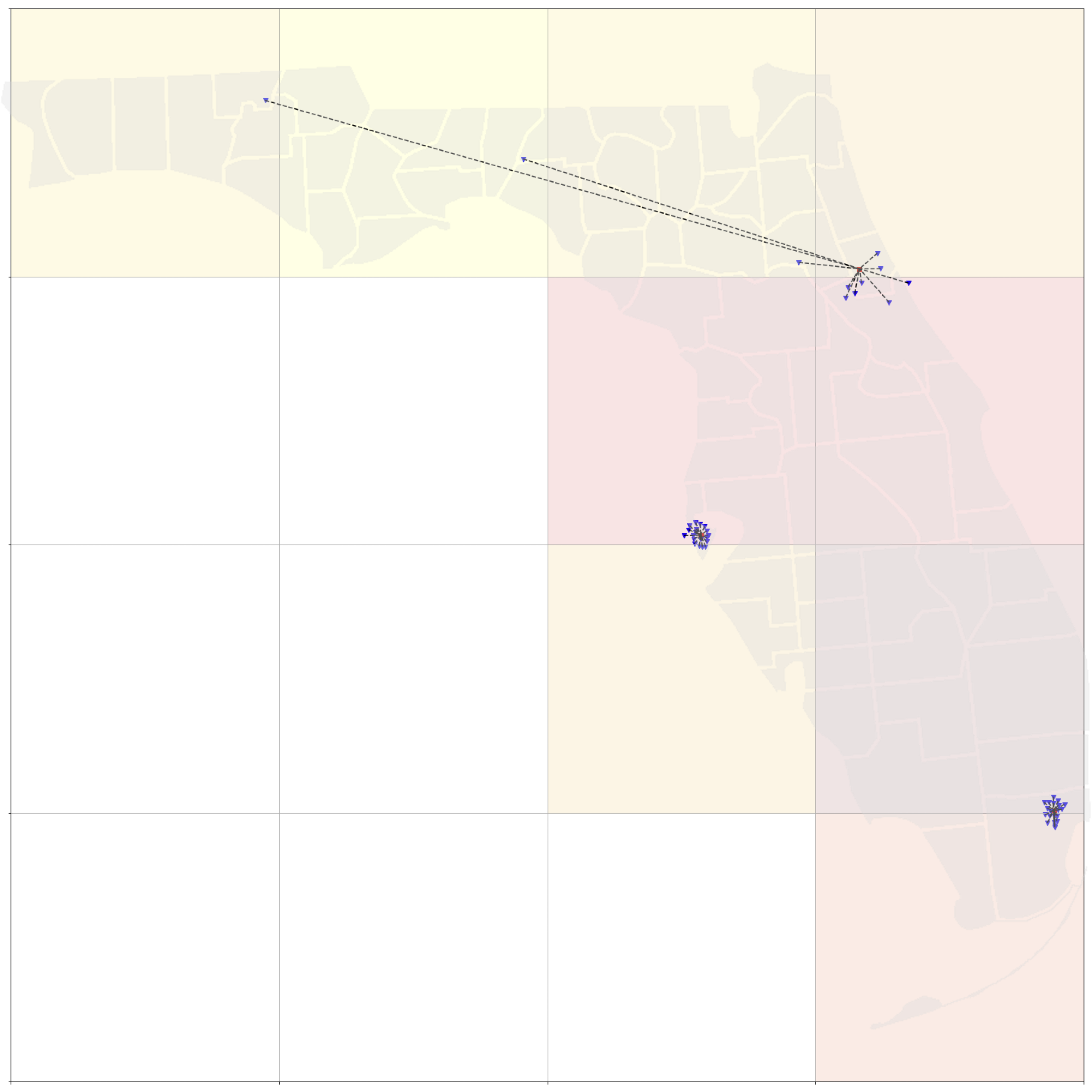}
     \end{subfigure}
     \hfill
     \begin{subfigure}[h]{0.45\textwidth}
         \centering
         \caption{Instance = ZIP, Method = Heuristic}
         \includegraphics[width=0.60\textwidth]{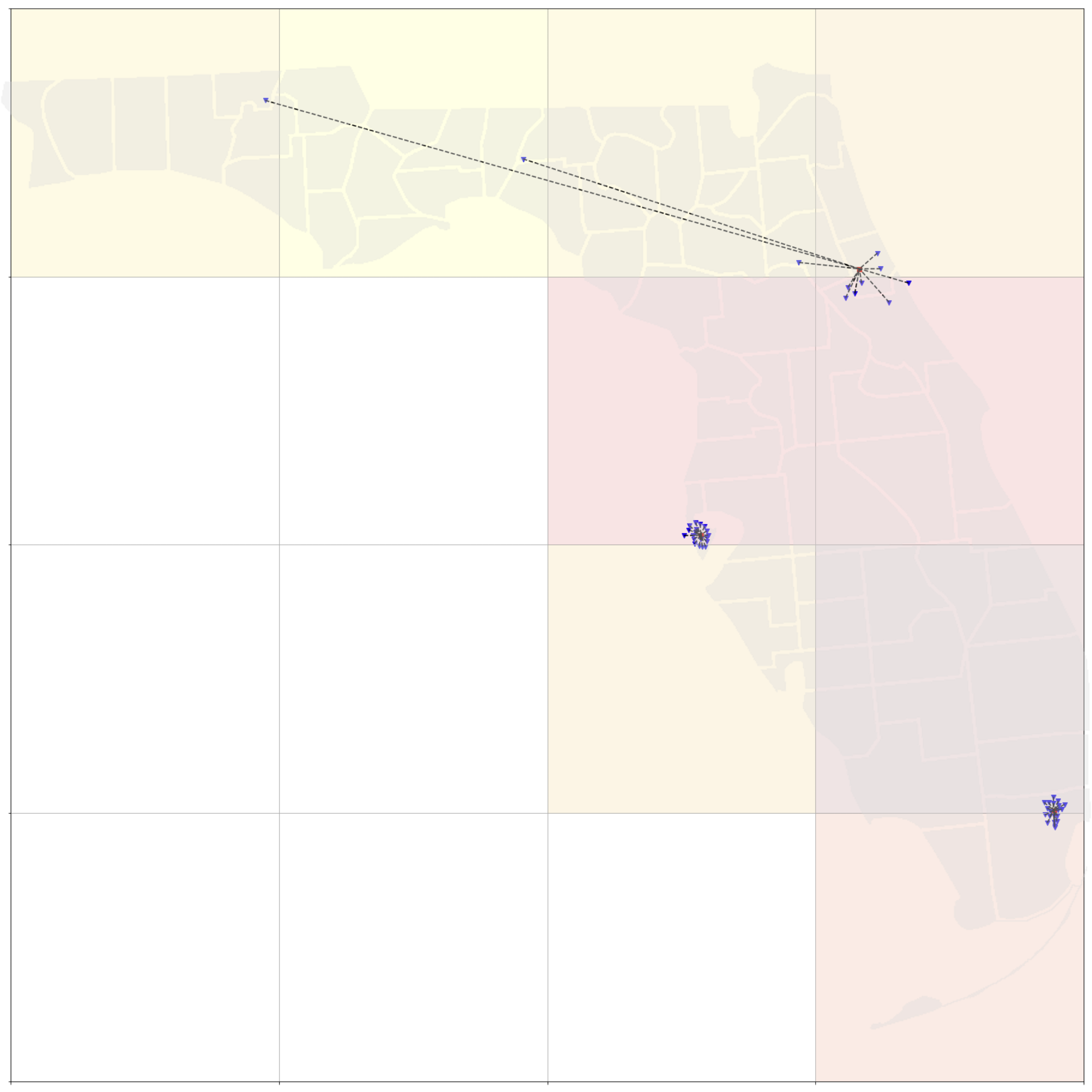}
     \end{subfigure}
\end{figure}

\FloatBarrier

\begin{figure}[h]
    \caption{Visualization of Solutions to the Small Instance and Grid Size 2 Under Center Objective}
    \label{Fig:CenterMap}
     \centering
     \begin{subfigure}[h]{0.45\textwidth}
         \centering
         \caption{Instance = City, Method = Exact}
         \includegraphics[width=0.60\textwidth]{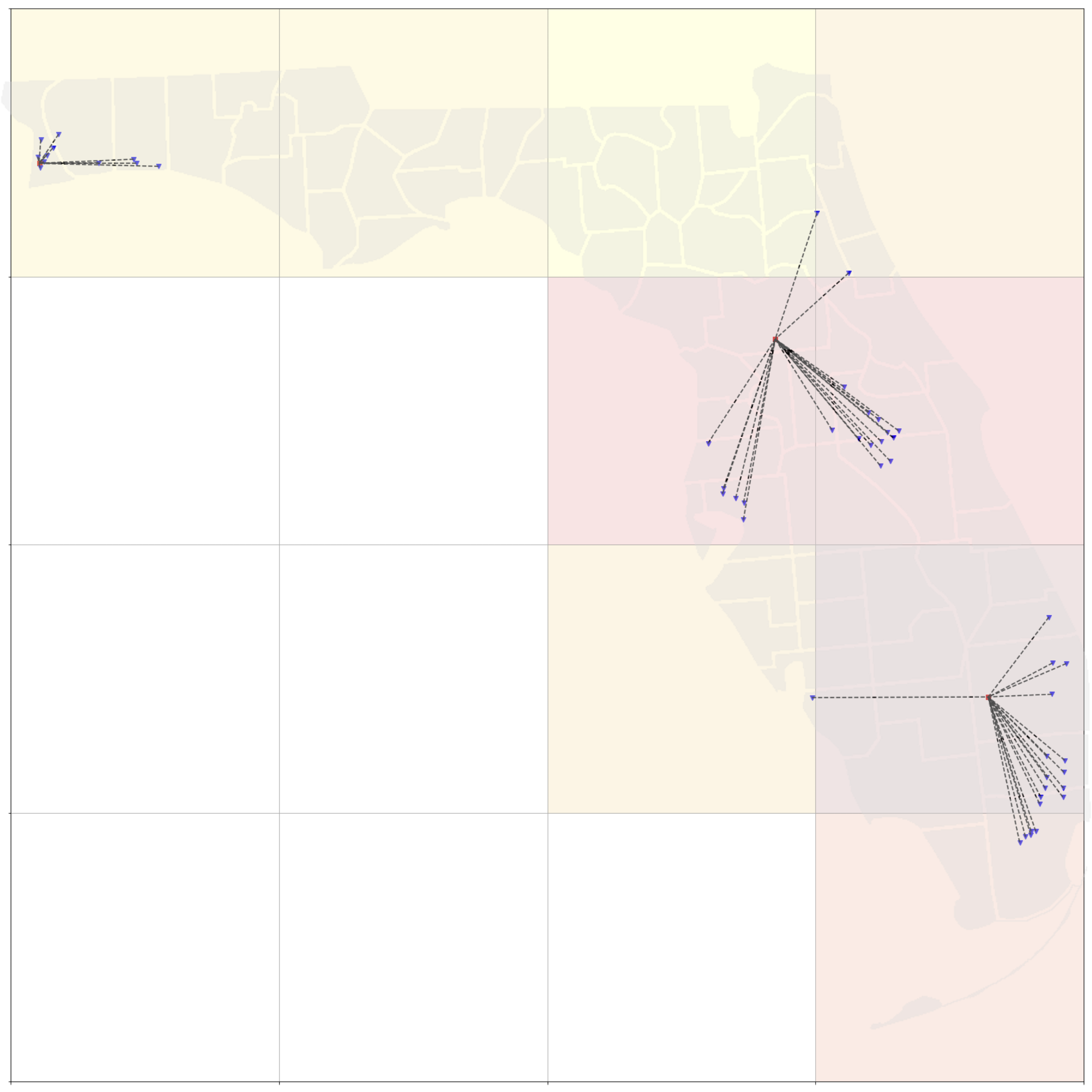}
     \end{subfigure}
     \hfill
     \begin{subfigure}[h]{0.45\textwidth}
         \centering
         \caption{Instance = City, Method = Heuristic}
         \includegraphics[width=0.60\textwidth]{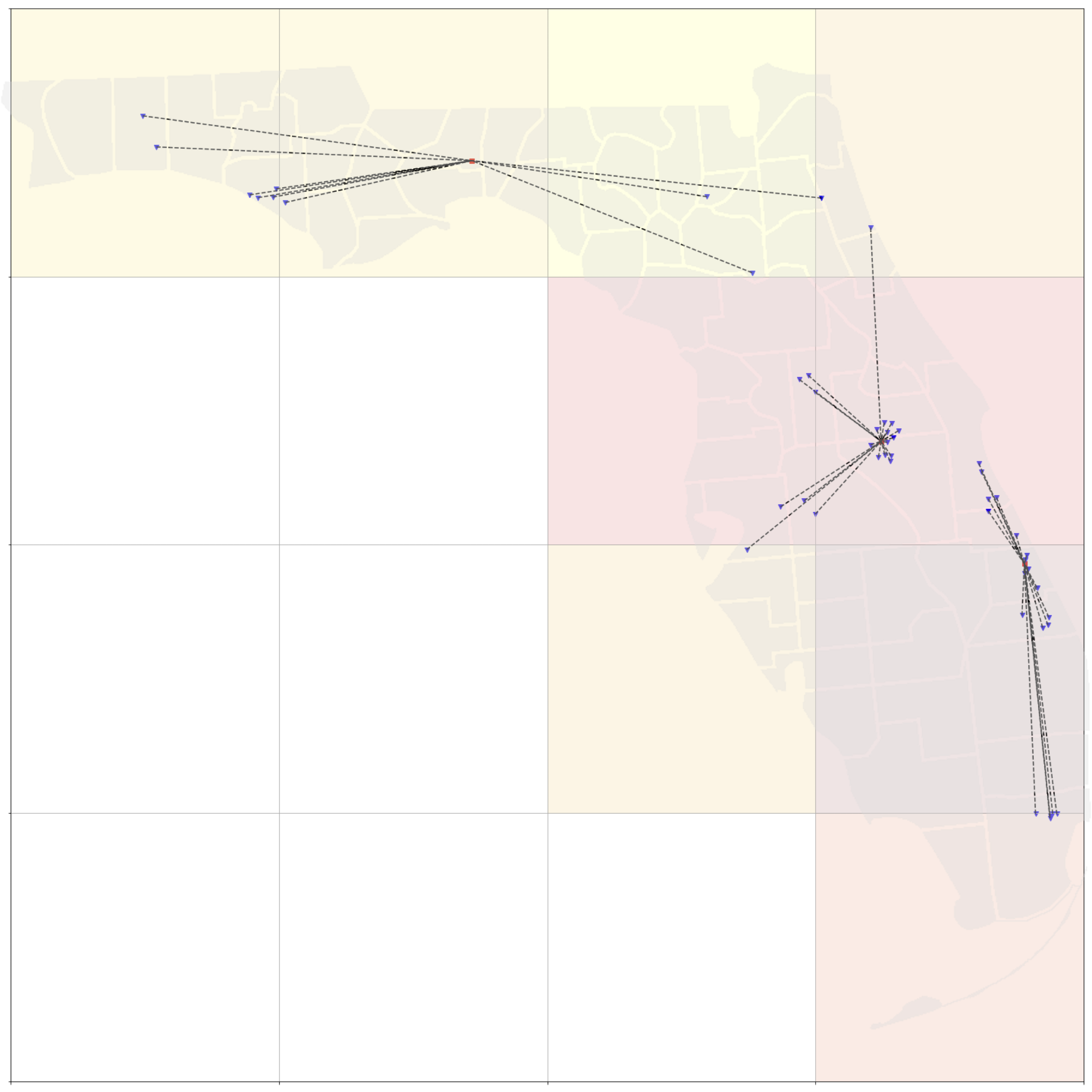}
     \end{subfigure} \\
     \centering
     \begin{subfigure}[h]{0.30\textwidth}
         \centering
         \caption{\tiny Instance = ZIP, Method = Exact}
         \includegraphics[width=0.9\textwidth]{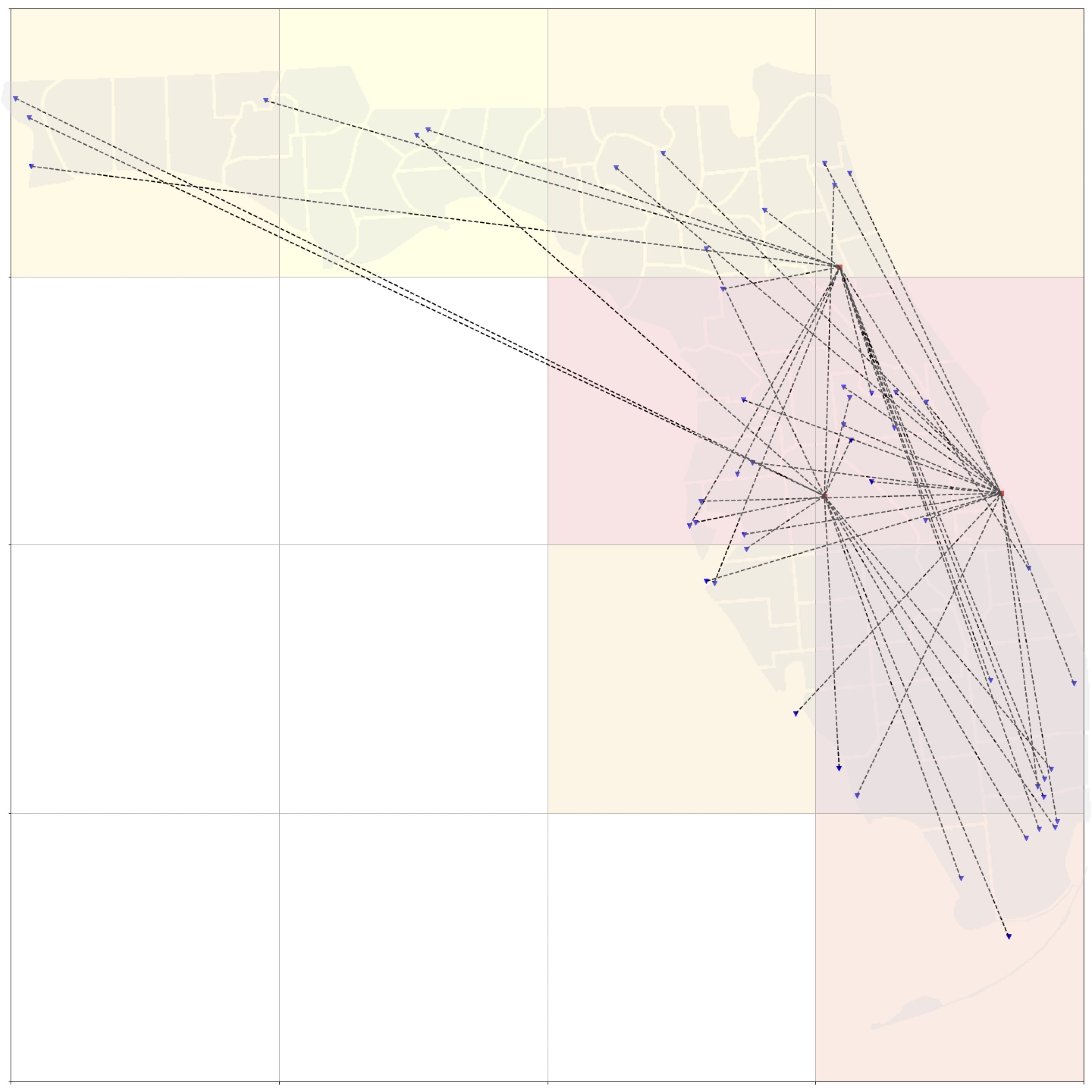}
     \end{subfigure}
     \hfill
     \begin{subfigure}[h]{0.30\textwidth}
         \centering
         \caption{\tiny Instance = ZIP, Method = Heuristic}
         \includegraphics[width=0.9\textwidth]{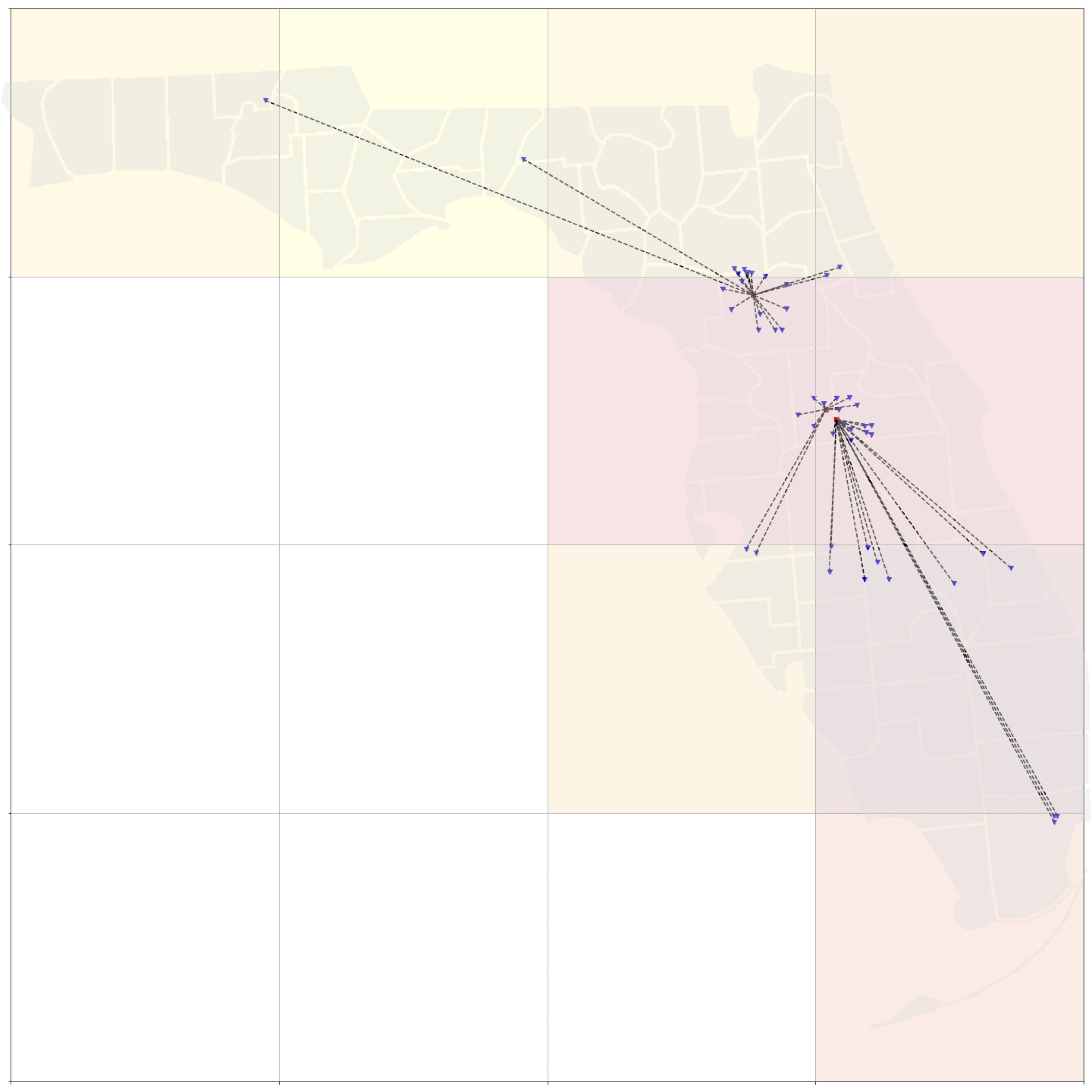}
     \end{subfigure}
     \hfill
     \begin{subfigure}[h]{0.30\textwidth}
         \centering
         \caption{\tiny Instance = ZIP, Method = Improvement}
         \includegraphics[width=0.9\textwidth]{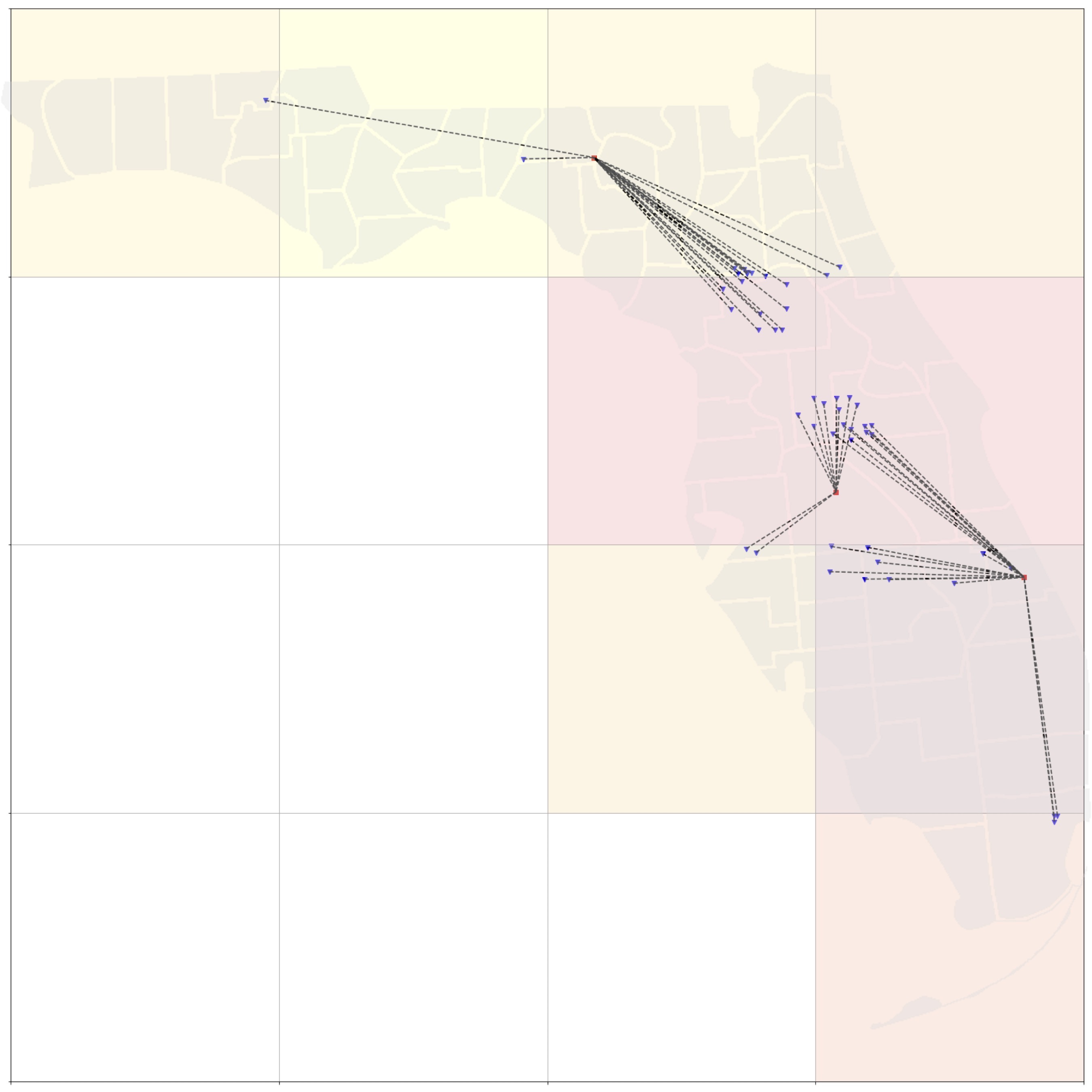}
     \end{subfigure}
\end{figure}

\FloatBarrier

Finally, we visualize the solutions on maps (Figures \ref{Fig:MedianMap} and \ref{Fig:CenterMap}), both for the small instance and with grid coefficient $g = 2$. The shading represents population percentiles (10, 30, 50, 70, 90\%), from yellow to red. Figure \ref{Fig:MedianMap} shows the solution map under the median objective. In each of the four cases, the exact and heuristic methods yield the same solution. Figure \ref{Fig:CenterMap} shows the solution map under the center objective. In the city case, the heuristic solution is $58\%$ worse than the exact solution; the difficulty comes from assigning demand in the northern regions. For the ZIP code instance, the exact method yields a highly suboptimal solution, while the heuristic yields a solution that appears reasonable. However, from the zone perspective, there are some situations resembling Figure \ref{Img:ImprovementCases}. After the improvement heuristic is applied, the objective improves by $35\%$. 

\FloatBarrier

\section{Conclusion}\label{sec:conclusion}
This paper proposes a novel extension for facility location problems by introducing the concept of demand location for binary demands. To best of our knowledge, this is the first work that formally conceptualizes joint facility and demand location. In addition, we propose a local-search network flow heuristic for the $k$-median and $k$-center problems with demand location, and evaluate its performance in a hurricane case study. The proposed local-search network flow framework is competitive with the exact solver in certain problem instances. Directions for future work include:
\begin{itemize}
    \item \textbf{Algorithm development:} There are several possible extensions of the network flow heuristic, such as random initialization, using different scoring rules, or adding additional types of local steps. It is also possible that heuristics developed for traditional facility location problems can be adapted to the joint location setting. We anticipate that this joint location problem will be a fertile ground for algorithm development.
    \item \textbf{Extensions and variations:} The current formulation assumes the demand nodes are binary and the distance is the only considered cost. The problem can be extended to cases where there is non-unit demand, location costs, or other objective functions.
    \item \textbf{Applications to real-life problems:} It would be interesting to further study the JFDLP in real-life contexts (i.e., refugee placement, sheltering against war, evacuation planning, homeless shelters). What is the gain of joint optimization compared to optimizing demand and facility nodes separately?
\end{itemize}

\section*{Acknowledgments}\label{sec:acknowledgments}

We sincerely thank Dr. Kelsey Rydland (Northwestern University Libraries) for providing the GIS data to assist this study, Dr. Samir Khuller (Northwestern Computer Science) for initial discussions about network flow heuristics, and Dr. Sanjay Mehrotra (Northwestern Industrial Engineering and Management Sciences) for funding A.K.K.

\newpage

\bibliographystyle{chicago}

\bibliography{references}

\end{document}